\documentclass[12pt,oneside]{amsart}
\usepackage{amsthm,amsfonts,amssymb,amscd,euscript}
\makeatletter
 \@addtoreset{equation}{subsection}
 \makeatother

\newtheorem{theorem}[subsection]{Theorem}
\newtheorem{conjecture}[subsection]{Conjecture}
\newtheorem{proposition}[subsection]{Proposition}
\newtheorem{proposition1}[subsection]{Proposition}
\newtheorem{lemma}[subsection]{Lemma}
\newtheorem{corollary}[subsection]{Corollary}
\newcommand{\para}{\refstepcounter{subsubsection}\par\smallskip\noindent
{\rm (\thesubsubsection)\hspace{0.2cm}}}
\newcommand{\parag}{\refstepcounter{subsection}\par\medskip\noindent
{\rm (\thesubsection)\hspace{0.2cm}}}

\theoremstyle{definition}
\newtheorem{definition}[subsection]{Definition}
\newtheorem{example}[subsection]{Example}
\newtheorem{examples}[subsection]{Examples}
\theoremstyle{remark}

\newtheorem*{remark*}{Remark}

\newcommand{\eref}[1]{{\rm{(\ref{#1})}}}
\newcommand{\mt}[1]{\operatorname{#1}}
\newcommand{\NE}{{\overline{NE}}}
\newcommand{\Diff}[2]{\operatorname{Diff}_{#1}{\left(#2\right)}}
\newcommand{\Supp}[1]{\operatorname{Supp}{\left(#1\right)}}
\newcommand{\Sing}{\operatorname{Sing}}
\newcommand{\Pic}{\operatorname{Pic}}
\newcommand{\Const}{\operatorname{Const}}
\newcommand{\compl}[1]{\operatorname{compl}{\left(#1\right)}}
\newcommand{\compll}[1]{\operatorname{compl}'{\left(#1\right)}}
\newcommand{\pair}[1]{\left(#1\right)}
\newcommand{\mult}{\operatorname{mult}}
\newcommand{\discr}[1]{\operatorname{discr}{\left(#1\right)}}
\newcommand{\totaldiscr}[1]{\operatorname{totaldiscr}{\left(#1\right)}}
\newcommand{\dis}[1]{a{\left(#1\right)}}

\newcommand{\pal}{\text{---}}

\newcommand{\fr}[1]{\left\{ #1 \right\}}
\newcommand{\ov}[1]{\overline{#1}}
\newcommand{\down}[1]{\left\lfloor #1\right\rfloor}
\newcommand{\up}[1]{\left\lceil #1\right\rceil}
\newcommand{\bir}{\dasharrow}

\newcommand{\CC}{{\mathbb C}}
\newcommand{\RR}{{\mathbb R}}
\newcommand{\ZZ}{{\mathbb Z}}
\newcommand{\QQ}{{\mathbb Q}}
\newcommand{\PP}{{\mathbb P}}
\newcommand{\NN}{{\mathbb N}}

\newcommand{\EE}{{\mathbb E}}
\newcommand{\OOO}{{\mathcal O}}
\newcommand{\KKK}{{\EuScript{K}}}
\newcommand{\NNN}{{\EuScript{N}}}
\newcommand{\MMM}{{\Phi}}
\newcommand{\Msm}{{\MMM}_{\mt{\mathbf{sm}}}}
\newcommand{\Mm}{{\MMM}_{\mt{\mathbf{m}}}}
\newcommand{\ep}{\varepsilon}
\newcommand{\var}{\varphi}

\newcommand{\De}{\Delta}

\title{The first main theorem on complements:
from global to local}
\author{Yu.~G.~Prokhorov}\thanks{The first author was
partially supported by the grant INTAS-OPEN-97-2072. The second
author was partially supported by the grant NSF-9800807.}
\address{Yu. Prokhorov: Department of Mathematics,
Tokyo Institute of Technology, Oh-Okayama, Meguro, Tokyo, Japan \&
Department of Mathematics (Algebra Section), Moscow State
University, 117234 Moscow, Russia}
\email{prokhoro@mech.math.msu.su}
\author{V.~V.~Shokurov}
\address{V.~V.~Shokurov: The Johns Hopkins University,
Department of Mathematics, Baltimore, Maryland,21218, USA}
\email{shokurov@math.jhu.edu}

\begin{document}

\begin{abstract}
The aim of this paper is to clarify and generalize techniques of
\cite{Sh1} (see also \cite{Pr1} and \cite{Pr2}). Roughly
speaking, we prove that for local Fano contractions the existence
of complements can be reduced to the existence of complements for
lower dimensional projective Fano varieties.
\end{abstract}

\maketitle
\section*{Introduction}
The aim of this paper is to clarify and generalize techniques of
\cite[Sect.~7]{Sh1} (see also \cite{Pr1}, \cite{Pr2}). We prove
that for local Fano contractions the existence of complements can
be reduced to the existence of complements for lower dimensional
projective Fano varieties. The main conjecture on $n$-complements
(Conjecture~\cite[1.3]{Sh1}) states that they are bounded in each
given dimension.
\par
Roughly speaking, an $n$-complement is a ``good'' member of the
multiple anti-log canonical linear system. A multitude of examples
support the conjecture \cite{Abe}, \cite{Is}, \cite{IP},
\cite{KeM}, \cite{Ko-SGT}, \cite{MP}, \cite{Pr2}, \cite{Sh},
\cite{Sh1}. As was noticed in \cite{Sh}, complements have good
structures which are related to restrictions of linear systems and
Kawamata-Viehweg vanishing. The latter essentially explains a
tricky structure of $n$-complement boundaries (cf. inequality in
\eref{def-coplements-n} below). In the main conjecture we consider
log pairs $(X/Z,D)$ consisting of Fano contractions $X/Z$ and
boundaries $D$. To use an induction in a proof of the conjecture
we need to divide log pairs and their complements into two types
with respect to the dimension of the base $Z$, namely, local
whenever $\dim(Z)>0$, and global otherwise. Equivalently, in the
global case $Z$ is a point and $X$ is a projective log Fano. We
prove, for local log Fano contractions, the existence of an
$n$-complement, where $n\in\NNN$ and the set $\NNN$ comes from
lower dimensional projective log Fano varieties. This is called
the first main theorem on complements (see
Theorem~\ref{result-Fano-1} below): from global to local. The
proof uses the LogMMP, so it is conditional in dimensions
$n=\dim(X)\ge 4$ and the proof for $n\le 3$. The core idea is to
extend an $n$-complement from a central fiber of a good
modification for $(X/Z,D)$ (cf. the proof of Theorem~5.6 and
Example~5.2 in \cite{Sh}). Moreover, such an approach allows us to
control some numerical invariants of complements: e.g.,
indices and their type, exceptional or non-exceptional, and their
regularity (cf. \cite[Sect. 7]{Sh1}).

The second theorem, from local to global, will be discussed in the
next paper. Its prototype is the global case in \cite{Sh1} (cf.
also tigers in \cite{KeM}) that uses local and inductive
complements \cite[Sect. 2]{Sh1}. An elementary but really generic
case of the second theorem is Theorem~\ref{result-Fano-2}. It is a
modification of the first one. Other cases show that
the main
difficulty of the Borisov-Alekseev conjecture (see
\ref{conjecture-boundedness-log-Fano}) concerns $\ep_d$-log
terminal log Fano varieties of dimension $d$, namely, that they
are bounded for some $\ep_d>0$ depending on the dimension $d$. For
instance, in the dimension $2$, $\ep_2=6/7$.

The paper is organized as follows. Section~1 is auxiliary. In
Section~2 we introduce the very important notion of exceptional
pairs. In Section~3 we prove the main result
(Theorem~\ref{result-Fano-1}). Some corollaries and applications
are discussed in Section~4. Finally, in Section~5 we present the
global version of Theorem~\ref{result-Fano-1}.

\section{Preliminaries}
\subsection*{Notation}\quad
\newline\noindent
\begin{tabular}{lp{12cm}}
 $\KKK(X)$& the function field of the variety $X$;\\
 $D_1\approx D_2$& prime divisors $D_1$, $D_2$ give the
 same discrete valuation of $\KKK(X)$;\\
 $K_X$& canonical (Weil) divisor, we will
frequently write $K$ if no confusion is likely. \\
\end{tabular}
\par\smallskip\noindent
All varieties are assumed to be algebraic and defined over $\CC$,
the field of complex numbers. A \textit{contraction} (or
\textit{extraction}, if we start with $X$ instead of $Y$) is a
projective morphism of normal varieties $f\colon Y\to X$ such that
$f_*\OOO_Y=\OOO_X$. A \textit{blow-up} is a birational extraction.
We will use the standard abbreviations and notation of Minimal
Model Program as MMP, lc, klt, plt, $\equiv$, $\sim$,
$\down{\cdot}$, $\up{\cdot}$, $\fr{\cdot}$, $\NE(X/Z)$,
$\dis{E,D}$, $\discr{X,D}$, $\totaldiscr{X,D}$; see \cite{KMM},
\cite{Ut}, \cite{Ko}. Everywhere below, if we do not specify the
opposite, a \textit{boundary} means a $\QQ$-boundary, i.e. a
$\QQ$-Weil divisor $D=\sum d_iD_i$ such that $0\le d_i\le 1$ for
all $i$. A \textit{log variety} (\textit{log pair}) $(X/Z\ni o,D)$
is, by definition, a contraction $X\to Z$ which is considered
locally near the fiber over $o\in Z$ and a boundary $D$ on $X$. By
the dimension of a log pair $(X/Z\ni o,D)$ we mean the dimension
of the total space $X$.
\par

\begin{definition}[\cite{Sh}]
\label{definition-complements}
Let $(X/Z,D)$ be a log variety. Then
\para
\label{def-coplements-num}
\textit{numerical complement} is an $\RR$-boundary $D'\ge D$, such
that $K+D'$ is lc and numerically trivial;
\para
\label{def-coplements-R}
$\RR$-\textit{complement} is an $\RR$-boundary $D'\ge D$ such that
$K+D'$ is lc and $\RR$-linearly trivial;
\para
\label{def-coplements-Q}
$\QQ$-\textit{complement} is a $\QQ$-boundary $D'\ge D$ such that
$K+D'$ is lc and $\QQ$-linearly trivial.
\para
\label{def-coplements-n}
Write $D=S+B$, where $S=\down{D}$, $B=\fr{D}$. Then an
$n$-\textit{complement} is a $\QQ$-boundary $D^+$ such that
$K+D^+$ is lc, $n(K+D^+)\sim 0$ and $nD^+\ge nS+\down{(n+1)D}$.
\end{definition}

Note that an $\RR$-complement can be considered as an
$n$-complement for $n=\infty$ because the limit of the inequality
in \eref{def-coplements-n} for $n\to \infty$ gives as $D'\ge D$.
All these definitions can be done in the more general situation:
when $D$ is an $\RR$-subboundary (i.e. an $\RR$-divisor $D=\sum
d_iD_i$ with $d_i\le 1$ for all $i$).

Obviously, there are the following implications:
\begin{quote}
$\exists$ $\QQ$-complement $\Longrightarrow$ $\exists$
$\RR$-complement $\Longrightarrow$ $\exists$ numerical complement.
\end{quote}
The simple example below shows that an $n$-complement is not
necessarily a $\QQ$-complement (even not a numerical complement).

\begin{example}
Let $P_1$, $P_2$, $P_3$ be a different points on $\PP^1$. Put
$D:=P_1+(\frac{1}{2}+\ep)P_2+(\frac{1}{2}-\ep)P_3$ and
$D':=P_1+\frac{1}{2}P_2+\frac{1}{2}P_3$, where $0<\ep\ll 1$. Then
$K+D'$ is a $2$-complement of the log divisor $K+D$. However
$D'\ge D$ is wrong, i.e. $K+D'$ is not a $\QQ$-complement of the
log divisor $K+D$.
\end{example}
Under additional restriction on coefficients of $D$ (for example,
if $D$ is \textit{standard}, see \eref{definition-NNN}) we have
$D^+\ge D$ in \eref{def-coplements-n}, see \cite[2.7]{Sh1} or
\cite{Pr3}. Therefore $D^+$ is a $\QQ$-complement in this case.

The question on the existence of complements naturally arises for
varieties of Fano or Calabi-Yau type, i.e. for varieties with nef
anti-log canonical divisor. However the nef property of $-(K+D)$
does not guarantee the existence of complements
\cite[1.1]{Sh1}.

\begin{proposition}[{\cite[5.5]{Sh}}]
Let $(X/Z\ni o,D)$ be a log variety. Assume that $K+D$ is lc and
$-(K+D)$ (semi)ample over $Z$. Then near the fiber over $o$ there
exists a $\QQ$-complement of the log divisor $K+D$.
\end{proposition}

\parag
\label{definition-NNN}
Fix a subset $\MMM\subset [0,1]$. We will write simply $D\in\MMM$
if all the coefficients of $D$ are contained in $\MMM$. For
example, we can consider $\MMM=\Msm:=\{1-1/m \mid
m\in\NN\cup\{\infty\}\}$ (this is called the case of
\textit{standard coefficients}). However some of our statements
and conjectures can be formulated for another choice of $\MMM$
(see \eref{eq-def-N-M} below).

\parag
\label{assumptions-NNN}
Let $(X,D)$ be a projective log variety such that:
\para
$K_X+D$ is lc;
\para
$D\in\MMM$;
\para
$-(K_X+D)$ nef and big;
\para
there exists some $\QQ$-complement of $K_X+D$ (this condition holds
if $-(K_X+D)$ is semi-ample, for example, by \cite[3-1-2]{KMM}
this holds if $K_X+D$ is klt).
\par\smallskip
Such a pair we call a \textit{log Fano variety}.

\parag
Notation as above. Define the
minimal complementary number by
\setcounter{equation}{\value{subsubsection}}
\begin{equation}
\label{eq-compl-def}
\compl{X,D}:=\min\{m\mid K_X+D\ \text{is $m$-complementary}\}.
\end{equation}
\setcounter{subsubsection}{\value{equation}}\noindent and consider
the set
\begin{multline*}
\NNN_d(\MMM):=\{m\in\NN\mid\exists\ \text{a log Fano variety} \
(X,D)\ \text{of dimension $d$}\\ \text{such that} \ D\in\MMM\
\text{and}\ \compl{X,D}=m\}.
\end{multline*}
For example, $\NNN_1([0,1])=\{1, 2, 3, 4, 6\}$ (see \cite{Sh}).
Taking products with $\PP^1$, one can show that
$\NNN_{d-1}(\MMM)\subset\NNN_d(\MMM)$. For inductive purposes we
put $\NNN_0([0,1])=\{1,2\}$.

By induction we define
\setcounter{equation}{\value{subsubsection}}
\begin{multline}
\label{eq-def-N-M}
\Mm^1:=\Msm,\qquad N_1:=\max\NNN_1(\Mm^1),\\
\Mm^d:=\Msm\cup\left[1-\frac{1}{N_{d-1}+1},1\right], \qquad
N_d=\max\left( \cup_{k=1}^d\NNN_d(\Mm^k) \right).
\end{multline}
\setcounter{subsubsection}{\value{equation}}\noindent We do not
exclude the case $N_d=\infty$ (and then $\Mm^d:=\Msm$), however,
we hope that $N_d<\infty$ (see
\ref{conjecture-boundedness-complements} below). By \cite[5.2]{Sh}
we have
\[
N_1=6,\qquad \Mm^2=\Msm\cup [6/7,1].
\]
It was proved in \cite{Sh1} that $N_2$ is finite. By construction,
$N_d\ge N_{d'}$ and $\Mm^d\subset\Mm^{d'}$ if $d\ge d'$.

\begin{lemma}[cf. {\cite[2.7]{Sh1}}]
\label{ge}
If $\alpha\in \Mm^d$, then for any $n\le N_{d-1}$ we have
\[
\down{(n+1)\alpha}\ge n\alpha.
\]
\end{lemma}
\begin{proof}
If $\alpha\in\Msm$, then $\alpha=1-1/m$ for some $m\in\NN$. In
this case we write $n\alpha=q+k/m$, where $q=\down{n\alpha}$ and
$k/m=\fr{n\alpha}$, $k\in\ZZ$, $0\le k\le m-1$. Then
\[
\down{(n+1)\alpha}=\down{q+k/m+1-1/m}=\left\{
\begin{array}{ll}
 q\qquad&\text{if}\ k=0,\\
 q+1\qquad&\text{otherwise.}
\end{array}
\right.
\]
In both cases $\down{(n+1)\alpha}\ge q+k/m=n\alpha$. Assume that
$\alpha\notin\Msm$. Then $\alpha>1-\frac{1}{N_{d-1}+1}$ and
\[
\down{(n+1)\alpha}\ge \down{n+1-\frac{n+1}{N_{d-1}+1}}\ge n\ge
n\alpha.
\]
\end{proof}

\begin{corollary}
\label{ge-1}
Let $(X,D)$ be a log pair such that $D\in \Mm^d$ and let $D^+$ be
an $n$-complement with $n\le N_{d-1}$. Then $D^+\ge D$.
\end{corollary}

\begin{lemma}[cf. {\cite[Lemma~4.2]{Sh}}]
\label{coeff-diff}
Let $(X,D)$ be a lc log pair, let $S:=\down{D}$ and $B:=\fr{D}$.
Assume that $K+S$ is plt and $D\in \Mm^{d}$ for some $d$ (resp.
$D\in \Msm$). Then $\Diff{S}{B}\in \Mm^{d}$ (resp. $\Diff{S}{B}\in
\Msm$).
\end{lemma}
\begin{proof}
Write $B=\sum b_jB_j$, $0<b_j<1$. Let $\alpha$ be a coefficient
$\Diff{S}{B}$. Then by \cite[3.10]{Sh},
\setcounter{equation}{\value{subsubsection}}
\begin{equation}
\label{eq-alpha}
\alpha=\frac{m-1}{m}+\sum_j\frac{b_jn_{j}}{m},
\end{equation}
\setcounter{subsubsection}{\value{equation}}\noindent where
$m\in\NN$, $n_{j}\in\NN\cup\{0\}$. Since $K_{S}+\Diff{S}{B}$ is lc
(see \cite[17.7]{Ut}), $\alpha\le 1$ and we may assume that
$\alpha<1$. Using $b_j\ge 1/2$ one can easily show that in
\eref{eq-alpha} $\sum n_j\le 1$ (see \cite[Lemma~4.2]{Sh}). If
$n_j=0$ for all $j$ in \eref{eq-alpha}, then, obviously,
$\alpha\in\Msm$. Otherwise $n_{j_0}=1$ for some $j_0$ and $n_j=0$
for $j\neq j_0$ in \eref{eq-alpha}. Then
$\alpha=\frac{m-1+b_{j_0}}{m}$. If $b_{j_0}\in\Msm$, then
$b_{j_0}=1-1/n$, $n\in\NN$ and $\alpha=\frac{mn-1}{mn}\in\Msm$. If
$b_{j_0}\ge 1-\frac{1}{N_{d-1}+1}$, then $\alpha\ge b_{j_0}\ge
1-\frac{1}{N_{d-1}+1}$. In both cases $\alpha\in\Mm^d$.
\end{proof}

\begin{conjecture}
\label{conjecture-boundedness-complements}
Notation as in \eref{assumptions-NNN}. Then $\NNN_d(\MMM)$ is
finite.
\end{conjecture}

The proof of Conjecture~\ref{conjecture-boundedness-complements}
in dimension two given in \cite{Sh1} relies heavily upon boundedness
results for log del Pezzo surfaces \cite{A}, see also \cite{N3}.
In arbitrary dimension there is the following

\begin{conjecture}
\label{conjecture-boundedness-log-Fano}
Fix $\ep>0$. Let $(X,D)$ be a normal projective log variety such
that:
\para
$K+D$ is $\QQ$-Cartier;
\para
$\totaldiscr{X,D}>-1+\ep$;
\para
$-(K_X+D)$ is nef and big.
\par\smallskip
Then $(X,\Supp{D})$ belongs to a finite number of algebraic
families.
\end{conjecture}
This conjecture is known to be true for $\dim(X)=2$. For
$\dim(X)\ge 3$ there are only particular results in this direction
\cite{B}, \cite{BB}. A new approach to the proof of
\ref{conjecture-boundedness-log-Fano} was proposed in \cite[Sect.
9]{KeM}.

\begin{conjecture}[Inductive Conjecture]
\label{conjecture-inductive-complements}
Let $(X,D)$ be such as in \eref{assumptions-NNN} (in particular
$D\in\MMM$). Assume that there exists a $\QQ$-complement of $K+D$
which is not klt. Then $K+D$ has an $n$-complement for
$n\in\NNN_{d-1}(\MMM)$. Moreover, this new complement also can be
taken non-klt.
\end{conjecture}
We may expect Conjecture~\ref{conjecture-inductive-complements}
for $\MMM=\Msm$ or $\MMM=\Mm^d$, where $d=\dim(X)$. In general, it
fails \cite[2.4]{Sh1}, \cite[8.1.2]{Pr3}. At the moment, this
conjecture is proved for $\dim(X)=2$ and $\MMM=\Mm^2$, \cite{Sh1}
(even in a stronger form).

\section{Exceptionality}

\begin{definition}
We say that a contraction $f\colon X\to Z$ is of \textit{local
type}, if $\dim(Z)>0$. Otherwise (i.e. $Z$ is a point) we say
that the contraction $f\colon X\to Z$ is of \textit{global type}.
\end{definition}
Thus a contraction of local type can be either birational or of
fiber type. In this case we are interested in the structure of
$f\colon X\to Z$ near the fixed fiber $f^{-1}(o)$, $o\in Z$ and
usually we assume that $X$ is a sufficiently small neighborhood of
the fiber over $o$.

\begin{definition}[{\cite[Sect. 5]{Sh}}, {\cite[1.5]{Sh1}}]
\label{def-exc-loc}
Let $(X/Z\ni o,\De)$ be a log variety of local type. Assume that
$K+\De$ has at least one $\QQ$-complement near the fiber over $o$.
Then $(X/Z\ni o,\De)$ is said to be \textit{exceptional} if for
any $\QQ$-complement $K+\De^+$ of $K+\De$ near the fiber over $o$
there exists at most one (prime) divisor $E$ of $\KKK(X)$ with
$\dis{E,\De^+}=-1$.
\end{definition}
Clearly, to be exceptional depends on the choice of
the base point $o\in Z$. As an immediate consequence of the
definition we have

\begin{lemma}
\label{lemma-crepant-ex}
Let $(X/Z\ni o,\De)$ and $(X'/Z\ni o,\De')$ be log varieties (of
local or global type) and let $f\colon X\to X'$ be a contraction
over $Z$. Assume that $K_{X'}+\De'$ is $\QQ$-Cartier and $\De$ is
a crepant pull back of $\De'$ (i.e. $f^*(K_{X'}+\De')=K_{X}+\De$
and $f_*\De=\De'$). Then $(X/Z\ni o,\De)$ is exceptional if and
only if $(X'/Z\ni o,\De')$ is.
\end{lemma}
\begin{proof}
Follows by \cite[3.10]{Ko}.
\end{proof}

\begin{proposition}
\label{two-divisors}
Let $(X/Z\ni o,\De)$ be a non-exceptional log variety of local
type and let $D$, $D'$ be $\QQ$-complements such that both $K+D$
and $K+D'$ are not klt. Let $S$ and $S'$ are divisors of $\KKK(X)$
such that $\dis{S,D}=-1$ and $\dis{S',D'}=-1$. Assume that
$S\not\approx S'$. Then there exists a $\QQ$-complement $G$ of
$K+\De$ such that $\dis{S,G}=\dis{E,G}=-1$ for some $E\not\approx
S$.
\end{proposition}
\begin{proof}[Proof (cf. {\cite[2.7]{MP}}, {\cite[2.4]{IP}})]
Note that $D'-D$ is $\QQ$-Cartier and numerically trivial over
$Z$. Put $D(\alpha):=D+\alpha(D'-D)$. Then $D(0)=D$, $D(1)=D'$ and
$K+D(\alpha)$ is a $\QQ$-complement for all $0\le\alpha\le 1$ (by
convexity of the lc property see \cite[1.4.1]{Sh} or
\cite[2.17.1]{Ut}). Fix an effective Cartier divisor $L$ on $Z$
(passing through $o$) and put $F:=f^*L$. For $0\le\alpha\le 1$,
define a function
\[
\varsigma(\alpha):=\sup\{\beta\mid K+D(\alpha)+\beta F\quad
\text{is lc}\},
\]
and put $T(\alpha):=D(\alpha)+\varsigma(\alpha)F$. Fix some log
resolution of $(X,D+D'+F)$ and let $\sum E_i$ be the union of the
exceptional divisor and the proper transform of $\Supp{D+D'+F}$.
Then $\varsigma(\alpha)$ can be computed as
\[
\varsigma(\alpha)=\max_{E_i}\{\beta\mid \dis{E_i,D(\alpha)+\beta
F}\ge -1\}.
\]
(see e.~g. \cite[0-2-12]{KMM}). In particular,
$\varsigma(\alpha)\in\QQ$. Hence $K+T(\alpha)$ is a
$\QQ$-complement. By the above, $\beta=\varsigma(\alpha)$ can be
computed from linear inequalities $\dis{E_i, D(\alpha)+\beta F}\ge
-1$, where $E_i$ runs through a finite number of prime divisors
$E_i$. Therefore the function $\varsigma(\alpha)$ is piecewise
linear and continuous in $\alpha$ and so are the coefficients of
$T(\alpha)$. By construction, $K+T(\alpha)$ is not klt for all
$0\le\alpha\le 1$. We claim that $\dis{S,T(0)}=-1$. Indeed,
$T(0)=D+\varsigma(0)F\ge D$. Thus $\dis{S,T(0)}\le \dis{S,D}=-1$.
Since $K+T(0)$ is lc, $\dis{S,T(0)}=-1$. Now, take
\[
\alpha_0:=\sup\{\alpha\mid \dis{S,T(\alpha)}=-1\}.
\]
By the above discussions $\alpha_0$ is rational (and
$\dis{S,T(\alpha_0)}=-1$). If $\alpha_0=1$, then we put $G:=T(1)$
and $E=S'$. Otherwise, for any $\alpha>\alpha_0$,
$\dis{S,T(\alpha)}>-1$. Hence there is a divisor $E\not\approx S$
of $\KKK(X)$ such that $\dis{E,T(\alpha)}=-1$. Again we can take
$E$ to be a component of $\sum E_i$. Thus $E$ does not depend on
$\alpha$ if $0<\alpha-\alpha_0\ll 1$. Obviously,
$\dis{E,T(\alpha_0)}=-1$ and we can put $G:=T(\alpha_0)$.
\end{proof}

\begin{corollary}
\label{two-divisors-non-except}
Let $(X/Z\ni o,\De)$ be a non-exceptional log variety of local
type, let $D\ge \De$ be a $\QQ$-complement such that $K+D$ is not
klt and let $S$ be a divisor of $\KKK(X)$ such that
$\dis{S,D}=-1$. Then there is a $\QQ$-complement $G\ge \De$ such
that $\dis{S,G}=\dis{E,G}=-1$ for some divisor $E\not\approx S$ of
$\KKK(X)$.
\end{corollary}
\begin{proof}
Since $(X/Z\ni o,\De)$ is non-exceptional, there is a
$\QQ$-complement $D'\ge \De$ such that $\dis{S',D'}=-1$ for some
$S'\not\approx S$. Then one can apply
Proposition~\ref{two-divisors}.
\end{proof}

\begin{corollary}
\label{two-divisors-except}
Let $(X/Z\ni o,\De)$ be an exceptional log variety of local type.
Then there exists a uniquely defined divisor $S$ of $\KKK(X)$
such that for any $\QQ$-complement $D$ one has $\dis{E,D}>-1$
whenever $E\not\approx S$ in $\KKK(X)$.
\end{corollary}

We call the divisor $S$ defined in \ref{two-divisors-except} the
\textit{central divisor} of an exceptional log pair $(X/Z\ni
o,\De)$.

\begin{corollary}
\label{dimension}
Let $(X/Z\ni o,\De)$ be a exceptional log variety of local type,
let $S$ be the central divisor. Then the center of $S$ on $X$ is
contained in the fiber over $o$.
\end{corollary}
\begin{proof}
Let $K+D$ be a $\QQ$-complement such that $\dis{S,D}=-1$ and let
$H$ be a general hyperplane section of $Z$ passing through $o$.
Since $f^*H$ does not contain the center of $S$, $\mult_Sf^*H=0$
and $\dis{S,D}=\dis{S,D+cf^*H}=-1$ for all $c$. Take $c$ so that
$K+D+cf^*H$ is maximally lc. Then, as in the proof of
Proposition~\ref{two-divisors}, $\dis{E,D+cf^*H}=-1$ for some
$E\not\approx S$, a contradiction.
\end{proof}

\begin{example}
\label{ex-sing}
Consider a log canonical singularity $X\ni o$ (i.e. $X=Z$ and
$\De=0$). Then it is exceptional if and only if for any boundary
$B$ on $X$ such that $K+B$ is lc there exists at most one divisor
$E$ of $\KKK(X)$ with $\dis{E,B}=-1$. For example, a
two-dimensional log terminal singularity is exceptional if and
only if it is of type $\EE_6$, $\EE_7$ or $\EE_8$ (see
\cite[5.2.3]{Sh}, \cite{MP}).
\end{example}

In the global case Definition~\ref{def-exc-loc} has a different
form:

\begin{definition}
Let $(X,\De)$ be a log variety of global type. Assume that $K+\De$
has at least one $n$-complement. Then $(X,\De)$ is said to be
\textit{exceptional} if any $\QQ$-complement $K+\De^+$ of $K+\De$
is klt (i.e. $\dis{E,\De^+}>-1$ for any divisor $E$ of
$\KKK(X)$).
\end{definition}

\begin{examples}
\label{ex-global}
(i) Let $X=\PP^1$, $Z=\mt{pt}$ and let $\De=\sum_{i=1}^r
(1-1/m_i)P_i$, $m_i\in\NN$, where $P_1$,\dots, $P_r$ are different
points. The divisor $-(K+\De)$ is nef if and only if $\sum_{i=1}^r
(1-1/m_i)\le 2$. In this case, the collection $(m_1,\dots,m_r)$
gives us an exceptional pair if and only if it is (up to
permutations) one of the following:
\[
\begin{array}{lllllllllll}
 &&&E_6:&(2,3,3)&&E_7:&(2,3,4)&&E_8:&(2,3,5)\\
 &&&\widetilde E_6:&(3,3,3)&
 \quad&\widetilde E_7:&(2,4,4)&\quad&\widetilde E_8:&(2,3,6)\\
 &&&\widetilde D_4:&(2,2,2,2)&&&&&&\\
\end{array}
\]
\par
(ii) Let $X=\PP^d$, $Z=\mt{pt}$ and let $\De=\sum_{i=1}^{d+2}
(1-1/m_i)\De_i$, $m_i\in\NN$, where $\De_1$,\dots, $\De_{d+2}$ are
hyperplanes in $\PP^d$. The log divisor $-(K+\De)$ is nef if and
only if $\sum 1/m_i\le 1$. If $(X,\De)$ is exceptional, then
$-(K+\De_j+\sum_{i\ne j}(1-1/m_i)\De_i)$ is not nef for all $j$.
Hence $\sum_{i\ne j} 1/m_i>1$. In this situation it is easy to
prove the existence of a constant $\Const(d)$ such that $m_j\le
\Const(d)$ for all $j$ (cf. \cite[8.16]{Ko}). Therefore there are
only a finite number of possibilities for exceptional collections
$(m_1,\dots,m_{d+2})$.
\end{examples}

Examples above and many other facts (see \cite{Sh1}, \cite{MP},
\cite{IP}, \cite{Pr2}, \cite{Is}) show that in general we may
expect the following principle:
\begin{itemize}
\item
non-exceptional pairs have good properties of $|-m(K+D)|$ for some
small $m$;
\item
exceptional pairs can be classified.
\end{itemize}

\section{Fano contractions}
In this section we prove Theorem~\ref{result-Fano-1} below. The
two dimensional version of this result was proved by the second
author in \cite{Sh}. Later it was generalized in \cite{Sh1},
\cite{Pr2}.
\begin{theorem}[Local case]
\label{result-Fano-1}
Let $\MMM:=\Mm^d$ (or $\MMM:=\Msm$) and let $(X/Z\ni o,D)$ be a
$d$-dimensional log variety of local type such that
\para
$D\in\MMM$;
\para
$K+D$ is klt;
\para
$-(K+D)$ is nef and big over $Z$.
\par \smallskip
Let $f\colon X\to Z$ be the structure morphism. Assume LogMMP in
dimension $d$. Then there exists a non-klt $n$-complement of $K+D$
near $f^{-1}(o)$ for $n\in\NNN_{d-1}(\MMM)$. Moreover, if $(X/Z\ni
o,D)$ is non-exceptional and
Conjecture~\ref{conjecture-inductive-complements} holds in
dimensions $d'\le d-1$ for $\MMM=\Mm^{d}$ (resp. $\MMM=\Msm$),
then $K+D$ is $n$-complementary near $f^{-1}(o)$ for
$n\in\NNN_{d-2}(\MMM)$. This complement also can be taken
non-exceptional.
\end{theorem}

In the non-exceptional case we expect more precise results. In this
case the existence of complements should depend on the topological
structure of the essential exceptional divisor (see
\cite[Sect.~7]{Sh1}).

\begin{example}
\label{ex-RDP}
Let $(Z\ni o)$ be a two-dimensional DuVal (RDP) singularity, let
$D=0$, and let $f=\mt{id}$. There is a non-klt $n$-complement of
$K_Z$ for some $n\in\NNN_1(\Msm)=\{1,2,3,4,6\}$ (see
\cite[5.2.3]{Sh}). The singularity is non-exceptional if it is of
type $A_n$ or $D_n$. In these cases there is a non-klt
$n$-complement for $n\in\NNN_0(\Msm)=\{1,2\}$.
\end{example}

The rough idea of the proof is very easy: we construct some
special blow-up of $X$ with irreducible exceptional divisor $S$
(Proposition~\ref{constr-plt-blowup}) and then apply inductive
properties of complements (Proposition~\ref{prodolj}) to reduce
the problem to a low dimensional (but possibly projective) variety
$S$.

\begin{lemma}
\label{weak-log-Fano}
Let $(X/Z,D)$ be a log variety such that $K_X+D$ is klt and
$-(K_X+D)$ nef and big over $Z$. Then there exists an effective
$\QQ$-divisor $D^{\mho}$ such that $K_X+D+D^{\mho}$ is again klt
and $-(K_X+D+D^{\mho})$ ample over $Z$.
\end{lemma}
\begin{proof}
Follows by Kodaira's lemma (see e.~g. \cite[0-3-3, 0-3-4]{KMM}).
\end{proof}

\begin{corollary}
\label{weak-log-Fano-cor}
Notation and assumptions as in \ref{weak-log-Fano}. Then the Mori
cone $\NE(X/Z)$ is polyhedral and generated by contractible
extremal rational curves.
\end{corollary}

\begin{definition}
\label{def_plt_blowup}
Let $(X,\De)$ be a log pair and let $g\colon Y\to X$ be a blow-up
such that the exceptional locus of $g$ contains exactly one
irreducible divisor, say $S$. Assume that $K_Y+\De_Y+S$ is plt and
$-(K_Y+\De_Y+S)$ is $g$-ample. Then $g\colon (Y\supset S)\to X$ is
called a \textit{purely log terminal (plt) blow-up} of $(X,\De)$.
\end{definition}
\textsc{Warning:} In contrast with log terminal modifications
\cite[3.1]{Sh2} purely log terminal blow-ups are not log crepant.

\begin{remark*}
Let $(X\ni o,D)$ be an exceptional singularity. Then by
Corollary~\ref{two-divisors-except} there is at most one plt
blow-up (see \cite[Prop.~6]{Pr1}).
\end{remark*}

\begin{proposition}[{\cite{Pr1}}, {\cite{Pr3}}, cf. {\cite{Sh3}}]
\label{constr-plt-blowup}
Let $(X,\De+\De^0)$ be a log variety such that $X$ is
$\QQ$-factorial, $\De\ge 0$, $\De^0\ge 0$, $K+\De+\De^0$ is lc but
not plt and $K+\De$ is klt. (We do not claim that $\De$ and
$\De^0$ have no common components). Assume LogMMP in dimension
$\dim(X)$. Then there exists a plt blow-up $g\colon (Y\supset
S)\to X$ of $(X,\De)$ such that
\para
\label{inductive-blow-up-i}
$K_Y+\De_Y+S+\De^0_Y=g^{*}(K+\De+\De^0)$ is lc;
\para
\label{inductive-blow-up-ii}
$K_Y+\De_Y+S+(1-\ep)\De^0_Y$ is plt and anti-ample over $X$ for
any $\ep>0$;
\para
\label{inductive-blow-up-iii}
$Y$ is $\QQ$-factorial and $\rho (Y/X)=1$.
\end{proposition}

Such a blow-up we call an \textit{inductive blow-up} of
$K_X+\De+\De^0$. It is important to note that this definition
depends on $\De$ and $\De^0$, not just on $\De+\De^0$. Such
blow-ups are very useful in the theory of complements. In the
local case one can construct a boundary $\De^0$ as in
Proposition~\ref{constr-plt-blowup} just by taking the pull-back of
some $\QQ$-divisor on $Z$. In the global case the problem of finding
$\De^0$ is not so easy.

\begin{proof}
First take a log terminal modification $h\colon V\to X$ of $(X,
\De+\De^0)$ (see \cite{Sh}, \cite[17.10]{Ut}). Write
\[
h^*(K+\De+\De^0)=K_V+\De_V+\De^0_V+E,
\]
where $\De_V$ and $\De^0_V$ are proper transforms of $\De$ and
$\De^0$, respectively, and $E$ is exceptional. One can take $h$ so
that $E$ is reduced and $E\ne 0$ (see \cite[17.10]{Ut},
\cite[3.1]{Sh2}). We claim that $K_V+\De_V+E$ cannot be nef over
$X$. Indeed, write
\[
h^*(K+\De)=K_V+\De_V+\sum \alpha_iE_i,\quad \text{where}\quad
\alpha_i<1\quad \text{for all}\quad i.
\]
This give us $h^*\De^0=\De^0_V+\sum (1-\alpha_i)E_i$, so
\[
K_V+\De_V+E\equiv-\De^0_V\equiv \sum (1-\alpha_i)E_i\quad
\text{over}\quad X,
\]
where $\sum (1-\alpha_i)E_i$ is effective, exceptional and $\ne
0$. This divisor cannot be $h$-nef (see e.~g. \cite[1.1]{Sh}).
Now, run $(K_V+\De_V+E)$-MMP over $X$. At the last step we get a
birational contraction $g\colon Y\to X$ which satisfies
\eref{inductive-blow-up-i}--\eref{inductive-blow-up-iii}.
\end{proof}

\para
\label{start}
We prove Theorem~\ref{result-Fano-1} by induction on $d$. So
assume that \ref{result-Fano-1} holds if $\dim(X)<d$. To begin the
proof, replace $X$ with its $\QQ$-factorialization (see
\cite[6.11.1]{Ut}). This preserves all our assumptions. Next, take
$D^{\mho}$ as in Lemma~\ref{weak-log-Fano} and put
$D^{\triangledown}:=D^{\mho}+cf^*H$, where $H$ is an effective
Cartier divisor on $Z$ passing through $o$ and $c$ is the log
canonical threshold $c=c(X,D+D^{\mho},f^*H)$ (the maximal such
that $K+D+D^{\mho}+cf^*H$ is lc). Then
\para
\label{triangld-lc}
$K+D+D^{\triangledown}$ is anti-ample over $Z$, lc and not klt.
\par\smallskip
Note that $D$ and $D^{\triangledown}$ can have common components.
Now, we distinguish two cases:
\begin{enumerate}
\item[(A)]
$K+D+D^{\triangledown}$ is plt (and $\down{D+D^{\triangledown}}\ne
0$);
\item[(B)]
$K+D+D^{\triangledown}$ is not plt.
\end{enumerate}
\par
In case (B) we consider an inductive blow-up $g\colon \widehat
X\to X$ of $(X,D+D^{\triangledown})$. Let $S$ be the (irreducible)
exceptional divisor. By \cite[5.4]{Sh} (or \cite[19.2]{Ut}) it is
sufficient to prove the existence of required complements on
$\widehat X$. Write
\setcounter{equation}{\value{subsubsection}}
\begin{equation}
\label{eq-first-def-D}
\begin{array}{l}
 g^*(K+D+D^{\triangledown})=K_{\widehat X}+\De+S+\widehat D^{\triangledown},\\
 g^*(K+D)=K_{\widehat X}+\De+aS,
\end{array}
\end{equation}
\setcounter{subsubsection}{\value{equation}}\noindent where
$\widehat D^{\triangledown}$ and $\De$ are proper transforms of
$D^{\triangledown}$ and $D$, respectively, and $a<1$. Note that
$\De+aS$ is not necessarily a boundary.
\par
In case (A) we put $\widehat X=X$, $g:=\mt{id}$,
$S=\down{D+D^{\triangledown}}$. In this case $S$ is irreducible by
the Connectedness Lemma \cite[17.4]{Ut} and because $S$ is normal
\cite[17.5]{Ut}. Define $\De$ from $D=\De+aS$, where $0\le a<1$
and $S$ is not a component of $\De$, and put $\widehat
D^{\triangledown}:=D+D^{\triangledown}-S-\De$. In both cases we
have by \eref{triangld-lc} and \eref{eq-first-def-D} the following
(see \cite[3.10]{Ko}):
\para
\label{condition-KDS}
$K_{\widehat X}+\De+S+\widehat D^{\triangledown}$ is lc, not klt,
$K_{\widehat X}+\De+aS$ is klt and both $-(K_{\widehat
X}+\De+S+\widehat D^{\triangledown})$ and $-(K_{\widehat
X}+\De+aS)$ are nef and big over $Z$.

\begin{lemma}
\label{M-klt}
Notation as above. There exist $\delta_0>0$ and a boundary $M$ on
$\widehat X$ such that
\para
\label{lemma-i}
$\De+aS\le M\le \De+S+(1-\delta_0)\widehat D^{\triangledown}$;
\para
\label{lemma-ii}
$K+M$ is klt;
\para
\label{lemma-iii}
$-(K+M)$ is nef and big over $Z$.
\par\smallskip
In particular, the Mori cone $\NE(\widehat X/Z)$ is polyhedral.
\end{lemma}
\begin{proof}
By \eref{triangld-lc}, $K+D+(1-\delta_0) D^{\triangledown}$ is klt
and anti-ample over $Z$ for sufficiently small positive
$\delta_0$. Take $M$ as the crepant pull-back
\setcounter{equation}{\value{subsubsection}}
\begin{multline}
\label{eq-multline-1}
K_{\widehat X}+M= g^*(K+D+(1-\delta_0) D^{\triangledown})=\\
g^*(K+D)+(1-\delta_0)\bigl(g^*(K+D+D^{\triangledown})-
g^*(K+D)\bigr)=\\ K_{\widehat
X}+\De+aS+(1-\delta_0)\bigl((K_{\widehat X}+\De+S+\widehat
D^{\triangledown}) -(K_{\widehat X}+\De+aS)\bigr).
\end{multline}
\setcounter{subsubsection}{\value{equation}}\noindent In other
words,
\begin{multline*}
M=\De + aS+(1-\delta_0)(S+\widehat D^{\triangledown}-aS)=\\
\De+\bigl(1-\delta_0(1-a)\bigr)S+(1-\delta_0) \widehat
D^{\triangledown}.
\end{multline*}\noindent
From \eref{eq-multline-1} we obtain that $K+M$ is klt
\cite[3.10]{Ko}, anti-nef and anti-big over $Z$. \eref{lemma-i}
holds if $a\le 1-\delta_0(1-a)$, i.e. for $0<\delta_0\ll 1$.
\end{proof}

\para
\label{plt-and-anti-ample}
Further, take $0<\lambda \ll \delta_0$ and put
\[
\widehat D^{\lambda}:=(1-\lambda)\widehat D^{\triangledown}.
\]
We claim that the log divisor $K_{\widehat X}+\De+S+\widehat
D^{\lambda}$ is plt and anti-ample over $Z$.
\par\smallskip
Indeed, in case (B), since $\rho (\widehat X/X)=1$, curves in the
fibers of $g$ generate an extremal ray, say $R$. Then $R\cdot
(K_{\widehat X}+\De+S+\widehat D^{\triangledown})=0$ (and
$K_{\widehat X}+\De+S+\widehat D^{\triangledown} $ is strictly
negative on all extremal rays $\ne R$, see \eref{triangld-lc} and
\eref{eq-first-def-D}). Further, by \eref{eq-first-def-D}
$\widehat D^{\triangledown}\equiv -(1-a)S$ over $X$ and this
divisor is positive on $R$. Thus $K_{\widehat X}+\De+S+\widehat
D^{\lambda}$ is strictly negative on all extremal rays of
$\NE(\widehat X/Z)$ for sufficiently small positive $\lambda$. By
Kleiman criterion, it is anti-ample. Finally, $K_{\widehat
X}+\De+S+\widehat D^{\lambda}$ is plt because $\widehat
D^{\lambda}\le \widehat D^{\triangledown}$. In case (A), our claim
obviously follows by \eref{triangld-lc}.
\par
Note that $M\le \De+S+\widehat D^{\lambda}$ by \eref{lemma-i}.
\par
Fix some set $F_1,\dots,F_r$ of prime divisors on $\widehat X$.
For $n\gg 0$, take a general member $F\in |-n(K_{\widehat
X}+\De+S+\widehat D^{\lambda})-\sum F_i|$ and put $B:=\widehat
D^{\lambda}+\frac1n(F+\sum F_i)$. We can take $F_1,\dots, F_r$ and
$n$ so that
\para
$K+\De+S+B$ is plt;
\para
\label{generators}
components of $B$ generate $N^1(\widehat X/Z)$.
\par \smallskip
By construction, we have
\para
$K+\De+S+B\equiv 0$ over $Z$.
\par \smallskip
Take $\ep>0$ so that $K+\De+S+(1+\ep)B$ is plt (see
\cite[2.17]{Ut}) and $M\le \De+S+(1-\ep)B$ (i.e. $1-\delta_0\le
(1-\ep)(1-\lambda)$, see proof of Lemma~\ref{M-klt}). Run
$(K+\De+S+(1+\ep)B)$-MMP over $Z$:
\[
 \vcenter{\hbox{
 \begin{picture}(215,80)
 \put(65,70){\vector(-2,-1) {40}}
 \put(25,30){\vector(2,-1){40}}
 \put(110,70){\line(2,-1){10}}
 \put(125,62.5){\line(2,-1){10}}
 \put(140,55){\vector(2,-1){10}}
 \put(150,30){\vector(-2,-1) {40}}
 \put(90,75){\makebox(0,0){$\widehat X$}}
 \put(10,40){\makebox(0,0){$X$}}
 \put(170,40){\makebox(0,0){$\ov X$}}
 \put(90,3){\makebox(0,0){$Z$}}
 \put(40,66){\makebox(0,0)[c]{\scriptsize $g$}}
 \put(40,13){\makebox(0,0)[c]{\scriptsize $f$}}
 \put(135,14){\makebox(0,0)[c]{\scriptsize $q$}}
 \end{picture}
 }}
\]
We will use $\ov\square$ to denote the proper transform on $\ov X$
of a divisor $\square$ on $\widehat X$. For each extremal ray $R$
we have $R\cdot B<0$ and $R\cdot (K+\De+S)>0$. Therefore any
contraction is either flipping or divisorial and contracts a
component of $B$. In particular, any divisorial contraction does
not contract $S$. At the end we get the situation when
$(K+\De+S+(1+\ep)B)$ is nef over $Z$ (we do not exclude the case
$\ov X=Z$). Since $K+\De+S+B\equiv 0$, $-(K+\De+S)$ is also nef
over $Z$.

\begin{lemma}
\label{NE}
We can run $(K+\De+S+(1+\ep)B)$-MMP so that on each step there is
a boundary $M\le \De+S+(1-\ep)B$ such that $K+M$ is klt and
$-(K+M)$ is nef and big over $Z$.
\end{lemma}
\begin{proof}
By Lemma~\ref{M-klt} such a boundary exists on the first step. If
$K+\De+S+(1+\ep)B \equiv \ep B$ is not nef over $Z$, then
$-(K+\De+S+(1-\ep)B)\equiv \ep B$ is also not nef over $Z$. Put
\[
t_0:=\sup\{t\mid -(K+M+t(\De+S+(1-\ep)B-M))\quad\text{is nef}\}.
\]
By Lemma~\ref{weak-log-Fano-cor} this supremum is a maximum and
is achieved on some extremal ray. Hence $t_0$ is rational and
$0< t_0<1$. Consider the boundary $M^0:=M+t_0(\De+S+(1-\ep)B-M)$.
Then $-(K+M^0)$ is nef over $Z$ and $M^0\le \De+S+(1-\ep)B$. We
claim that $-(K+M^0)$ is also big over $Z$. Assume the opposite.
By the Base Point Free Theorem, $-(K+M^0)$ is semi-ample over $Z$ and
defines a contraction $\var\colon \widehat X\to W$ onto a
lower-dimensional variety. Let $C$ be a general curve in a fiber.
Then $C\cdot (K+M^0)=C\cdot (K+\De+S+B)=0$, so $C\cdot
(\De+S+B-M^0)=0$. Since $C$ is nef, $\ep C\cdot B\le C\cdot
(\De+S+B-M^0)=0$ and $C\cdot B=0$. By \eref{generators}, $C\equiv
0$, a contradiction.
\par
Further, $\NE(\widehat X/Z)$ is polyhedral, so there is an
extremal ray $R$ such that $R\cdot (K+M^0)=0$ and $\ep R\cdot
B=-R\cdot (K+\De+S+(1-\ep)B)<0$. Hence $R\cdot
(K+\De+S+(1+\ep)B)<0$. Let $h\colon \widehat X\to Y$ be the
contraction of $R$. Put $M^0_Y:=h_*M^0$. Then
$K+M^0=h^*(K_Y+M^0_Y)$. Therefore $K_Y+M^0_Y$ is $\QQ$-Cartier,
klt and $-(K_Y+M^0_Y)$ is nef and big over $Z$. If $g$ is
divisorial, we can continue the process replacing $\widehat X$
with $Y$ and $M$ with $M=M^0_Y$. Assume that $g$ is a flipping
contraction and let
\[
\vcenter{\hbox{
 \begin{picture}(85,50)
 \put(14,35){\vector(1,-1){21}}
 \put(66,35){\vector(-1,-1){21}}
 \put(15,40){\line(1,0){6}}
 \put(26,40){\line(1,0){6}}
 \put(37,40){\line(1,0){6}}
 \put(48,40){\line(1,0){6}}
 \put(59,40){\vector(1,0){6}}
 \put(7,40){\makebox(0,0){$\widehat X$}}
 \put(77,40){\makebox(0,0){$X^+$}}
 \put(40, 6){\makebox(0,0){$Y$}}
 \put(19,20){\makebox(0,0){\scriptsize$h$}}
 \put(62.1,20.1){\makebox(0,0){\scriptsize$h^+$}}
 \end{picture}}}
\]
be the flip. Take $M:={h^+}^{-1}(M^0_Y)$. Again we have that
$-(K_{X^+}+{M}^+)=-{h^+}^*(K_Y+M^0_Y)$ is nef and big over $Z$.
Thus we can continue the process replacing $X$ with $X^+$.
\end{proof}
Finally, we get on $\ov X$
\para
\label{ov-plt}
$K+\ov\De+\ov S$ is plt;
\para
\label{ov-nef}
$-(K+\ov\De+\ov S)$ is nef over $Z$.

\begin{lemma}
\label{ov-big}
Notation as above. Then $-(K+\ov\De+\ov S)$ is semi-ample over
$Z$. Moreover, if $-(K+\ov\De+\ov S)$ is not ample, then it
defines a birational contraction over $Z$ with the exceptional
locus contained in $\Supp{\ov B}$. In particular, $-(K_{\ov
S}+\Diff{\ov S}{\ov\De}) =-(K+\ov\De+\ov S)|_{\ov S}$ is big (and
nef) over $q(\ov S)$.
\end{lemma}
\begin{proof}
By Lemma~\ref{NE} and the Base Point Free Theorem, $-(K+\ov\De+\ov S)$
is semi-ample. Thus for some $n\in\NN$ the linear system
$|-n(K+\ov\De+\ov S)|$ defines a contraction $\ov X\to W$. For any
curve $C$ in a fiber we have $C\cdot \ov B=0$. Since the
components of $\ov B$ generate $N^1(\ov X/Z)$ (see
\eref{generators}), we have that $C\cdot \ov B_i<0$ for some
component $\ov B_i$ of $\ov B$. Hence $C\subset\Supp{\ov B}$.
\end{proof}
Note that $q\colon \ov S\to q(\ov S)$ is also a contraction:

\begin{lemma}
$q_*\OOO_{\ov S}=\OOO_{q(\ov S)}$ and $q(\ov S)=f(g(S))$ is
normal.
\end{lemma}
\begin{proof}
See the proof of Lemma~3.6 in \cite{Sh}.
\end{proof}

By Lemma~\ref{coeff-diff}, $\Diff{\ov S}{\ov\De}\in\MMM$ (recall
that we put $\MMM=\Mm^d$ or $\MMM=\Msm$).

\begin{lemma}
\label{return}
Assume that near $q^{-1}(o)$ there exists an $n$-complement
$K_{\ov S}+\Diff{\ov S}{\ov\De}^+$ of $K_{\ov S}+\Diff{\ov
S}{\ov\De}$. Then near $q^{-1}(o)$ there exists an $n$-complement
$K+D^+$ of $K+D$. Moreover, if $K_{\ov S}+\Diff{\ov S}{\ov\De}^+$
is not klt, then $K+D^+$ is not exceptional.
\end{lemma}
\begin{proof}
By Proposition~\ref{prodolj} any $n$-complement of $K_{\ov
S}+\Diff{\ov S}{\ov\De}$ can be extended to an $n$-complement of
$K+\ov\De+\ov S$. By \ref{bir-prop} we can pull-back complements
of $K+\De+S$ under divisorial contractions because they are
$(K+\De+S)$-positive. Finally, note that the proper transform of
an $n$-complement under a flip is again an $n$-complement. Indeed,
the inequality in \eref{def-coplements-n}, obviously, is preserved
under any birational map which is an isomorphism in codimension
one. The log canonical property is preserved by \cite[2.28]{Ut}.
\end{proof}

\begin{lemma}
\label{dim-q(S)>0}
If $\dim(q(\ov S))>0$, then
\para
$(X/Z\ni o,D)$ is not exceptional;
\para
there is a non-exceptional $n$-complement of $K+D$ with
$n\in\NNN_{d-2}(\MMM)$.
\end{lemma}
\begin{proof}
(i) follows by Corollary~\ref{dimension}. Note that $(\ov S/q(\ov
S)\ni o,\Diff{\ov S}{\ov\De})$ satisfies the conditions of our
theorem (see Lemma~\ref{coeff-diff}). By inductive hypothesis we
may assume that there is a non-klt $n$-complement of $K_{\ov
S}+\Diff{\ov S}{\ov\De}$ for $n\in\NNN_{d-2}(\MMM)$. The rest
follows by Lemma~\ref{return}.
\end{proof}

\parag
\label{non-exceptional-}
Going back to the proof of Theorem~\ref{result-Fano-1}, assume
that $(X/Z\ni o, D)$ is non-exceptional (i.e. there exists a
non-exceptional complement $K+D+\Upsilon$) and $q(\ov S)=o$. We
have to show only that there exists a non-exceptional
$n$-complement of $K+D$ with $n\in\NNN_{d-2}(\MMM)$. By
Lemma~\ref{dim-q(S)>0} we may assume that $q(\ov S)=o$, i.e. $\ov
S$ is projective. By Corollary~\ref{two-divisors-non-except} we
can take $\Upsilon$ so that $\dis{S,D+\Upsilon}=-1$ (and
$\dis{E,D+\Upsilon}=-1$ for some $E\not\approx S$). Let
$\widehat\Upsilon$ and $\ov\Upsilon$ be proper transforms of
$\Upsilon$ on $\widehat X$ and $\ov X$, respectively. Then
\[
g^*(K+D+\Upsilon)=K_{\widehat X}+\De+S+\widehat\Upsilon.
\]
Moreover, $\dis{E,\De+S+\widehat\Upsilon}=\dis{E,\ov\De+\ov
S+\ov\Upsilon}=-1$, because $K_{\widehat
X}+\De+S+\widehat\Upsilon\equiv 0$. Thus $K_{\ov X}+\ov\De+\ov
S+\ov\Upsilon$ is not plt (near $q^{-1}(o)$).

\begin{lemma}
\label{non-klt}
Assumptions as in \eref{non-exceptional-}. Then $K_{\ov
S}+\Diff{\ov S}{\ov\De+\ov\Upsilon}$ is not klt.
\end{lemma}
\begin{proof}
By the Adjunction \cite[17.6]{Ut} it is sufficient to prove that
$K+\ov\De+\ov S+\ov\Upsilon$ is not plt near $\ov S$. Taking into
account discussions above, we see that this is a consequence of
Lemma~\ref{connect} below.
\end{proof}

By Lemma~\ref{non-klt} and
Conjecture~\ref{conjecture-inductive-complements} we obtain that
there is a non-klt $n$-complement of $K_{\ov S}+\Diff{\ov
S}{\ov\De}$ with $n\in\NNN_{d-2}(\MMM)$. By Lemma~\ref{return}
this proves Theorem~\ref{result-Fano-1}.

The following example illustrates the proof of
Theorem~\ref{result-Fano-1}:
\begin{example}
As in Example~\ref{ex-RDP}, let $(Z\ni o)$ be a two-dimensional
DuVal (RDP) singularity, let $D=0$, and let $f=\mt{id}$. In this
case, $g\colon \widehat X\to X$ is a weighted blow-up (with
suitable weights) and $\widehat X\bir \ov X$ is the identity map.
Hence $S\simeq\PP^1$. Write $\Diff S{0}=\sum_{i=1}^r
(1-1/m_i)P_i$, where $P_1,\dots,P_r$ are different points. We have
the following correspondence between types of $(Z\ni o)$ and
collections $(m_1,\dots,m_r)$ (see \ref{ex-sing} and
\ref{ex-global}):
\par\medskip
\begin{center}
\begin{tabular}{c||c|c|c|c|c}
$(Z\ni o)$&$A_n$&$D_n$&$E_6$&$E_7$&$E_8$\\ \hline
$(m_1,\dots,m_r)$&$r\le
2$&$(2,2,m)$&$(2,3,3)$&$(2,3,4)$&$(2,3,5)$\\
\end{tabular}
\end{center}
\par\medskip\noindent
Thus $(Z\ni o)$ is exceptional if and only if it is of type $E_6$,
$E_7$ or $E_8$.
\end{example}

\begin{lemma}[see \cite{Pr2}, cf. {\cite[6.9]{Sh}},
{\cite[Proposition~2.1]{F}}]
\label{connect}
Let $(X/Z\ni o,D)$ be a be a log variety and let $f\colon X\to Z$
be the structure morphism. Assume that
\para
$K+D$ is lc and not plt near $f^{-1}(o)$;
\para
$K+D\equiv 0$ over $Z$;
\para
there is an irreducible component $S\subset\down{D}$ such that
$f(S)\ne Z$.
\par \smallskip
Assume also LogMMP in dimension $\dim(X)$. Then $K+D$ is not plt
near $S\cap f^{-1}(o)$.
\end{lemma}

\begin{corollary}
\label{main-corollary}
Notation as in Theorem~\ref{result-Fano-1}. The following are
equivalent:
\para
\label{(i)}
$(X/Z\ni o, D)$ is an exceptional pair (of local type);
\para
\label{(ii)}
$q(\ov S)=o$ and $(\ov S,\Diff{\ov S}{\ov D})$ is an exceptional
pair (of global type).
\end{corollary}
\begin{proof}
\eref{(i)} $\Longrightarrow$ \eref{(ii)} follows by
Lemma~\ref{dim-q(S)>0} and Lemma~\ref{return}. The inverse
implication follows by Lemma~\ref{non-klt}.
\end{proof}

Define
\[
\compll{X,D}:=\min\{m\mid \ \text{there is non-klt $m$-complement
of $K+D$}\}.
\]
\begin{corollary}
Notation and assumptions as in Theorem~\ref{result-Fano-1}. Assume
that $(X/Z\ni o, D)$ is exceptional. Then
\[
\compll{X,D}=\compl{\ov S,\Diff{\ov S}{\ov D}}.
\]
\end{corollary}
\begin{proof}
The inequality $\le$ follows by Lemma~\ref{return}, so we show
$\ge$. Let $K+D^+$ be a non-klt $n$-complement of $K+D$. Then
$D^+\ge D$. By Corollary~\ref{two-divisors-except}
$\dis{S,D^+}=-1$. Consider the crepant pull-back
$g^*(K+D^+)=K_{\widehat X}+\De+S+\Upsilon$ and let $\ov \Upsilon$
be the proper transform of $\Upsilon$ on $\ov X$. Then $K_{\ov
S}+\Diff{\ov S}{\ov \De+\ov \Upsilon}$ is an $n$-complement of
$K_{\ov S}+\Diff{\ov S}{\ov \De}$.
\end{proof}
Note that for non-exceptional contractions we have only
$\compll{X,D}\le \compl{\ov S,\Diff{\ov S}{\ov D}}$:

\begin{example}
Let $(X\ni o)$ be a terminal $cE_8$-singularity given by the
equation $x_1^2+x_2^3+x_3^5+x_4^r=0$, $\gcd(r,30)=1$ and let
$g\colon (\widehat X,S)\to X$ be the weighted blow-up with weights
$(15r,10r,6r,30)$. Then $S=\PP^2$ and $\Diff S0=\frac12 L_1+
\frac23 L_2+\frac45 L_3+\frac{r-1}r L_4$, where $L_1,\dots,L_4$
are lines on $\PP^2$ in general position. Then $\compll{X,0}=1$
because $(X\ni o)$ is a $cDV$-singularity. On the other hand,
$\compl{\ov S,\Diff{\ov S}0}=6$.
\end{example}

\section{Exceptional Fano contractions}
In this section we study exceptional Fano contractions such as in
Theorem~\ref{result-Fano-1}.
\begin{proposition}
\label{prop-exc-1}
Notation and assumptions as in Theorem~\ref{result-Fano-1}. Assume
also Conjecture~\ref{conjecture-boundedness-complements} in
dimensions $\le d-1$. Let $(X/Z\ni o, D)$ is exceptional. Then
\[
\dis{E,D}\ge -1+\delta_d\quad\text{for any}\quad E\not\approx S,
\]
where $\delta_d>0$ is a constant which depends only on $d$.
\end{proposition}
\begin{proof}
Let $K+D^+$ be a non-klt $n$-complement with
$n\in\NNN_{d-1}(\MMM)$. Then $D^+\ge D$ (see \ref{ge-1}). By
definition of exceptional contractions $\dis{S,D^+}=-1$ and
$\dis{E,D^+}>-1$ for all $E\not\approx S$. Hence
$\dis{E,D^+}\ge-1+1/n$ because $n\dis{E,D^+}$ is an integer. Since
$D^+\ge D$, $\dis{E,D}\ge \dis{E,D^+}$. Thus we can take
$\delta_d:=1/\max(\NNN_{d-1}(\Mm^{d-1}))$.
\end{proof}

Assuming $D\in\Msm$, we obtain
\begin{corollary}
Notation and assumptions as in \ref{prop-exc-1}. Let $D_i$ be a
component of $D$ and let $d_i=1-1/m_i$ be its coefficient. If
$D_i\not\approx S$, then $m_i\le 1/\delta(d)$ and therefore there
is a finite number of possibilities for $d_i$.
\end{corollary}

\begin{corollary}[cf. \cite{Ko-SGT}]
Assume $LogMMP$ in dimensions $\le d$ and
Conjecture~\ref{conjecture-inductive-complements} in dimensions
$\le d-1$. Let $(X\ni o)$ be a $d$-dimensional klt singularity and
let $F=\sum F_i$ be an effective reduced Weil $\QQ$-Cartier
divisor on $X$ passing through $o$. Then we have either
$c_o(X,F)=1$ or $c_o(X,F)\le 1-1/N_{d-1}$, where $c_o(X,F)$ is the
log canonical threshold of $(X,F)$ \cite{Sh} (see also \cite{Ko})
and $N_{d-1}$ is such as in \eref{eq-def-N-M}.
\end{corollary}
This corollary is non-trivial only if
Conjecture~\ref{conjecture-boundedness-complements} holds in
dimension $\le d-1$.
\begin{proof}
Put $c:=c_o(X,F)$ and assume that $1-1/N_{d-1}<c<1$.
Theorem~\ref{result-Fano-1} gives us that there is an
$n$-complement $K+B$ of $K+cF$, where $n\le N_{d-1}$. Let $c^+_i$
be the coefficient of $F_i$ in $B$. By \eref{def-coplements-n},
$c^+\ge 1$. Hence $F\le B$ and $K+F$ is lc, a contradiction.
\end{proof}
In the case $1-1/(N_{d-2}+1)\le c=c_o(X,F)<1$, the pair $(X,cF)$
is exceptional. We expect that there are only a finite number of
possibilities for $c\in\left[1-1/(N_{d-2}+1,1)\right]$ in any
dimension. For example, this method gives us (see e.~g.
\cite[6.1.3]{Pr3}) that in dimension $d=2$ the set of all values
of $c_o(X,F)$ in the interval $[2/3,1]$ is
$\left\{2/3,7/{10},3/4,5/6,1\right\}$.

\begin{theorem}
\label{theorem-cover}
Fix $\ep>0$. Let $(X/Z\ni o, D)$ be a $d$-dimensional log variety
of local type such that
\para
$D\in\Msm$ (i.e. $D=\sum (1-1/m_i)D_i$, where
$m_i\in\NN\cup\{\infty\}$ and $D_i$'s are prime divisors);
\para
\label{assumption-totaldiscr}
$\totaldiscr{X,D}>-1+\ep$;
\para
$-(K+D)$ is nef and big over $Z$;
\para
$(X/Z\ni o, D)$ is exceptional.
\par\smallskip
Let $\var\colon X'\to X$ be a finite cover such that
\para
$X'$ is normal and irreducible;
\para
$\var$ is \'etale in codimension one outside of $\Supp{D}$;
\para
\label{assumption-divides}
the ramification index of $\var$ at the generic point of
components of $\var^{-1}(D_i)$ divides $m_i$.
\par \smallskip
Assume also LogMMP in dimensions $\le d$ and
Conjectures~\ref{conjecture-boundedness-log-Fano},
\ref{conjecture-boundedness-complements} and
\ref{conjecture-inductive-complements} for $\Msm$ in dimension
$d-1$. Then the degree of $\var$ is bounded by a constant
$\Const(d,\ep)$.
\end{theorem}
\begin{proof}
We will use notation of the proof of Theorem~\ref{result-Fano-1}.
Taking the fiber product with $\QQ$-factorialization we can reduce
the situation to the case when $X$ is $\QQ$-factorial. Note also
that $\var^{-1}\circ f^{-1}(o)$ is connected (because $X$ is
considered as a germ near $f^{-1}(o)$ and $X'$ is irreducible).
Consider the base change
\[
\label{CD-main}
\begin{CD}
 X'@>\var>>X\\
 @Vf'VV @VfVV\\
 Z'@>\pi>>Z\\
\end{CD}
\]
where $X'\stackrel{f'}\to Z'\stackrel{\pi}\to Z$ is the Stein
factorization. Then $f'\colon X'\to Z'$ is a contraction and
$\pi\colon Z'\to Z$ is a finite morphism. Define $D'$ and
$D^{\triangledown\prime}$ by
\setcounter{equation}{\value{subsubsection}}
\begin{equation}
\label{eq-var-preimage}
\begin{array}{l}
K_{X'}+D'=\var^*(K+D),\\
K_{X'}+D'+D^{\triangledown\prime}=\var^*(K+D+D^{\triangledown})
\end{array}
\end{equation}
\setcounter{subsubsection}{\value{equation}}\noindent (see
\cite[Sect.~2]{Sh}). This means that, for example, the coefficient
of a component $D_{i,j}'$ of $\var^{-1}(D_i)$ in $D'$ is as
follows
\[
d'_{i,j}=1-r_{i,j}(1-(1-1/m_i)),
\]
where $r_{i,j}$ is the ramification index at the generic point of
$D_{i,j}'$. Then by \eref{assumption-divides}, $D'\in \Msm$ (and
$D^{\triangledown\prime}\ge 0$). Obviously,
$K_{X'}+D'+D^{\triangledown\prime}$ is ample over $Z'$.

\para
First we consider case (A) (i.e. when $K+D+D^{\triangledown}$ is
plt, $S:=\down{D+D^{\triangledown}}\ne 0$, $\widehat X=X$,
$g:=\mt{id}$). We put
$S':=\down{D'+D^{\triangledown\prime}}=\var^{-1}(S)$. By
Lemma~\ref{dim-q(S)>0}, $S$ is compact and $S\subset f^{-1}(o)$.
Applying \cite[Sect.~2]{Sh} (or \cite[20.3]{Ut}) we get that
$K_{X'}+D'+D^{\triangledown\prime}$ is plt. By the Connectedness Lemma
\cite[17.4]{Ut} and the Adjunction \cite[17.6]{Ut}, $S'$ is connected,
irreducible and normal. Define $\De'$ from $D'=\De'+a'S'$, where
$0\le a'<1$. Let $\ov{X}'$ be the normalization of $\ov{X}$ in the
function field of $X'$. There is the commutative diagram
\[
 \vcenter{\hbox{
 \begin{picture}(98,50)
 \put(28,5){\vector(1,0){52}}
 \put(88,37){\line(0,-1){4}}
 \put(88,30){\line(0,-1){4}}
 \put(88,23){\line(0,-1){4}}
 \put(88,16){\vector(0,-1){4}}
 \put(20,37){\line(0,-1){4}}
 \put(20,30){\line(0,-1){4}}
 \put(20,23){\line(0,-1){4}}
 \put(20,16){\vector(0,-1){4}}
 \put(28,45){\vector(1,0){52}}
 \put(20,45){\makebox(0,0){$X'$}}
 \put(88,45){\makebox(0,0){$X$}}
 \put(88,5){\makebox(0,0){$\ov X$}}
 \put(20,5){\makebox(0,0){$\ov X'$}}
 \put(54,10){\makebox(0,0){\scriptsize$\ov\var$}}
 \put(59,50){\makebox(0,0){\scriptsize$\var$}}
 \put(13,25){\makebox(0,0){\scriptsize$\psi '$}}
 \put(96,25){\makebox(0,0){\scriptsize$\psi$}}
 \end{picture}
 }}
 \]
where $\ov \var\colon \ov{X}'\to \ov{X}$ is a finite morphism and
$\psi\colon X\bir \ov X$ and $\psi'\colon X'\bir \ov X'$ are
birational maps such that both $\psi^{-1}$ and $\psi^{\prime -1}$
do not contract divisors. Hence $\ov\var$ has the ramification
divisor only over $\Supp{\psi_*(D)}\subset\ov S\cup\Supp{\ov \De}$
and the ramification index of $\ov\var$ at the generic point of a
component over $\psi_*(D_i)$ is equal the ramification index of
$\var$ at the generic point of the corresponding component over
$D_i$. Applying $\psi_*$ and $\psi_*'$ to \eref{eq-var-preimage},
we obtain
\setcounter{equation}{\value{subsubsection}}
\begin{equation}
\begin{array}{l}
\label{eq-gather-X-D}
 K_{\ov X'}+\ov D'=\ov \var^*\left(K_{\ov X}+\ov D\right),\\
 K_{\ov X'}+\ov D'+\ov D^{\triangledown\prime}=
\ov \var^*\left(K_{\ov X}+\ov D+\ov D^{\triangledown}\right),
\end{array}
\end{equation}
\setcounter{subsubsection}{\value{equation}}\noindent where $\ov
D':=\psi_*' D'$ and $\ov
D^{\triangledown\prime}:=\psi_*'D^{\triangledown\prime}$. Recall
that $\ov S=\down{\ov D+\ov D^{\triangledown}}$ is irreducible.
Now, \eref{eq-gather-X-D} yields
\setcounter{equation}{\value{subsubsection}}
\begin{equation}
\label{eq-X-prime}
K_{\ov X'}+\ov \De'+\ov S'=\ov \var^*\left(K_{\ov X}+\ov \De+\ov
S\right),
\end{equation}
\setcounter{subsubsection}{\value{equation}}\noindent where $\ov
\De':=\psi_*' \De'$ and $\ov S':=\psi_*' S'$. By
\cite[Sect.~2]{Sh} and \eref{ov-plt} (see also \cite[20.3]{Ut}),
$K_{\ov X'}+\ov \De'+\ov S'$ is plt. Moreover, \eref{ov-nef} and
Lemma~\ref{ov-big} give us that $-(K_{\ov X'}+\ov \De'+\ov S')$ is
nef and big over $Z'$. It is sufficient to prove the boundedness
of the degree of the restriction $\ov \phi=\ov\var|_{\ov S'}
\colon\ov S'\to \ov S$. Indeed, $\deg \var=(\deg \ov \phi)r$,
where $r$ is the ramification index over $S$. By
\eref{assumption-totaldiscr} and \eref{assumption-divides}, $r$ is
bounded. Now, we consider log pairs $\pair{\ov S,\Diff{\ov
S}{\ov\De}}$ and $\pair{\ov S',\Diff{\ov S'}{\ov\De'}}$.
\par
Restricting \eref{eq-X-prime} on $\ov S$, we obtain
\[
K_{\ov S}+\Diff{\ov S}{\ov\De}=\ov \phi^*\left(K_{\ov S'}+
\Diff{\ov S'}{\ov\De'}\right).
\]
In particular, $\pair{K_{\ov S}+\Diff{\ov
S}{\ov\De}}^{d-1}=\pair{\deg \ov \phi}\pair{K_{\ov S'}+ \Diff{\ov
S'}{\ov\De'}}^{d-1}$. Both sides of this equality are positive by
Lemma~\ref{ov-big}.
\para
\label{para-A}
By the proof of Theorem~\ref{result-Fano-1}, there is an
$n$-complement $K_{\ov X}+\ov \De+\ov S+\ov \Upsilon$ of $K_{\ov
X}+\ov \De+\ov S$ with $n\le \max\NNN_{d-1}(\Msm)<\infty$. Define
$\ov \Upsilon'$ from
\[
\label{eq-X-prime-Up}
K_{\ov X'}+\ov \De'+\ov S'+\ov \Upsilon'= \ov \var^*\left(K_{\ov
X}+\ov \De+\ov S+ \ov \Upsilon\right),
\]
and put $\Theta:=\Diff{\ov S}{\ov\De+\ov \Upsilon}$ and
$\Theta':=\Diff{\ov S'}{\ov\De'+\ov \Upsilon'}$. Then $K_{\ov
S}+\Theta$ and $K_{\ov S'}+\Theta'$ are $n$-complements. Since
$K_{\ov S}+\Theta$ is klt (see Corollary~\ref{main-corollary}), we
have
\[
\totaldiscr{\ov S, \Diff{\ov S}{\ov\De}}\ge -1+\frac1n\ge
-1+\beta, \quad \text{where} \quad
\beta=\frac{1}{\max\NNN_{d-1}(\Msm)}.
\]
Similarly,
\[
\totaldiscr{\ov S', \Diff{\ov S'}{\ov\De'}}\ge -1+\beta.
\]
By \ref{conjecture-boundedness-log-Fano}, $\left(\ov
S,\Supp{\Diff{\ov S}{\ov\De}}\right)$ and $\left(\ov
S',\Supp{\Diff{\ov S'}{\ov\De'}}\right)$ belongs to a finite
number of algebraic families. Taking into account that $\Diff{\ov
S}{\ov\De}, \Diff{\ov S'}{\ov\De'}\in\Msm$ (see
Lemma~\ref{coeff-diff}) and $\Diff{\ov S}{\ov\De}\le \Theta$,
$\Diff{\ov S'}{\ov\De'}\le \Theta'$, we see that so are $\pair{\ov
S,\Diff{\ov S}{\ov\De}}$ and $\pair{\ov S',\Diff{\ov
S'}{\ov\De'}}$. This gives us that $\deg\ov \phi$ is bounded.

Now, we consider case (B). Let $\widehat X'$ be the normalization
of a dominant component of $\widehat X\times_X X'$ and let $S'$ be
the proper transform of $S$ on $\widehat X'$. We claim that
$g'\colon (\widehat X'\supset S')\to X'$ is a plt blow-up of
$(X',D')$. Consider the base change
\setcounter{equation}{\value{subsubsection}}
\begin{equation}
\label{CD1}
\begin{CD}
 \widehat X'@>\widehat \var>>\widehat X\\
 @Vg'VV @VgVV\\
 X'@>\var>>\phantom{.}X.\\
\end{CD}
\end{equation}
\setcounter{subsubsection}{\value{equation}}\noindent It is clear
that $\widehat\var\colon \widehat X'\to \widehat X$ is finite and
its ramification divisor can be supported only in $S\cup
\Supp{D}$. Then $S'$ is the exceptional divisor of the blow-up
$g'\colon \widehat X'\to X'$. We have
\setcounter{equation}{\value{subsubsection}}
\begin{equation}
\label{for-2}
K_{\widehat X'}+\De'+S'=\widehat \var^*(K_{\widehat X}+\De+S),
\end{equation}
\setcounter{subsubsection}{\value{equation}}\noindent where $\De'$
is a boundary. This divisor is plt \cite[2.2]{Sh}, \cite[20.3]{Ut}
and anti-ample over $X'$. By the Adjunction \cite[17.6]{Ut}, $S'$
is normal. On the other hand, $S'$ is connected near the fiber
over $o'\in Z'$. Indeed, $-\pair{K_{\widehat X'}+\widehat
D'+\widehat D^{\triangledown\prime}}$ is nef and big over $Z'$, by
\eref{CD1} and \eref{condition-KDS}. Since
$S'\subset\down{\widehat D'+\widehat D^{\triangledown\prime}}$, it
is connected by the Connectedness Lemma \cite[17.4]{Ut}. This
proves our claim.
\par
Now, as in case (A) we consider the commutative diagram
\[
 \vcenter{\hbox{
 \begin{picture}(98,50)
 \put(28,5){\vector(1,0){52}}
 \put(88,37){\line(0,-1){4}}
 \put(88,30){\line(0,-1){4}}
 \put(88,23){\line(0,-1){4}}
 \put(88,16){\vector(0,-1){4}}
 \put(20,37){\line(0,-1){4}}
 \put(20,30){\line(0,-1){4}}
 \put(20,23){\line(0,-1){4}}
 \put(20,16){\vector(0,-1){4}}
 \put(28,45){\vector(1,0){52}}
 \put(20,45){\makebox(0,0){$\widehat X'$}}
 \put(88,45){\makebox(0,0){$\widehat X$}}
 \put(88,5){\makebox(0,0){\phantom{.}$\ov X$.}}
 \put(20,5){\makebox(0,0){$\ov X'$}}
 \put(54,10){\makebox(0,0){\scriptsize$\ov\var$}}
 \put(59,50){\makebox(0,0){\scriptsize$\widehat\var$}}
 \put(13,25){\makebox(0,0){\scriptsize$\psi '$}}
 \put(96,25){\makebox(0,0){\scriptsize$\psi$}}
 \end{picture}
 }}
\]
Similar to case (A), $\pair{\ov S,\Diff{\ov S}{\ov\De}}$ and
$\pair{\ov S',\Diff{\ov S'}{\ov\De'}}$ are bounded. Hence we may
assume that $\deg \ov \phi$ is bounded, where $\ov
\phi=\ov\var|_{\ov S'} \colon\ov S'\to \ov S$. It remains to show
that the ramification index $r$ of $\ov\var$ at the generic point
of $S'$ is bounded. Clearly, $r$ is equal to the ramification
index of $\widehat\var$ at the generic point of $\widehat S'$.
Similar to \eref{eq-first-def-D} write
\setcounter{equation}{\value{subsubsection}}
\begin{equation}
\label{eq-a-prime}
\begin{array}{l}
 g^*(K_{X'}+D')=K_{\widehat X'}+\De'+a'S'.
\end{array}
\end{equation}
\setcounter{subsubsection}{\value{equation}}\noindent Then
\setcounter{equation}{\value{subsubsection}}
\begin{equation}
\label{eq-a-final}
1-a'=r(1-a)\ge r(1+\discr{X,D})> r\ep
\end{equation}
\setcounter{subsubsection}{\value{equation}}\noindent (see
\cite[Sect.~2]{Sh} or \cite[proof of 20.3]{Ut}). We claim that
$(S', \Diff{S'}{\De'})$ belong to a finite number of algebraic
families. Note that we cannot apply
\ref{conjecture-boundedness-log-Fano} directly because
$-(K_{S'}+\Diff{S'}{\De'})$ is not necessarily nef. As in case
(A), take $n$-complement $K_{\widehat X}+\De+S+\widehat\Upsilon$
with $n\le \max\NNN_{d-1}(\Msm)$. Similar to
\eref{eq-var-preimage} define $\widehat\Upsilon'$ and $\widehat
D^{\lambda\prime}$ (see \eref{plt-and-anti-ample}):
\begin{gather*}
K_{\widehat X'}+\De'+S'+\widehat\Upsilon'=
\widehat\var^*(K_{\widehat X}+\De+S+\widehat\Upsilon)\\
K_{\widehat X'}+\De'+S'+\widehat D^{\lambda\prime}=
\widehat\var^*(K_{\widehat X}+\De+S+\widehat D^{\lambda})\  .
\end{gather*}
Then $K_{S'}+\Diff{S'}{\De'+\widehat\Upsilon'}\equiv 0$ and by
\eref{plt-and-anti-ample}, $K_{S'}+\Diff{S'}{\De'+\widehat
D^{\lambda\prime}}$ is anti-ample. Hence
$-\pair{K_{S'}+\Diff{S'}{\De'+\alpha \widehat D^{\lambda\prime}+
(1-\alpha)\widehat\Upsilon'}}$ is ample for any $\alpha>0$. Note
that
\[
\totaldiscr{S', \Diff{S'}{\De'+\widehat\Upsilon'}}\ge -1+1/n.
\]
Thus we can apply Conjecture~\ref{conjecture-boundedness-log-Fano}
to $\pair{S', \Diff{S'}{\De'+\alpha \widehat D^{\lambda\prime}+
(1-\alpha)\widehat\Upsilon'}}$ for small positive $\alpha$. We
obtain that $S'$ is bounded. Now, as in \eref{para-A} we see that
so is $(S', \Diff{S'}{\De'})$. Take a sufficiently general curve
$\ell$ in a general fiber of $g'|_{S'}\colon S'\to g'(S')$. From
\eref{eq-a-prime} we have
\setcounter{equation}{\value{subsubsection}}
\begin{equation}
\label{eq-diff-ell-de}
-(K_{S'}+\Diff{S'}{\De'})\cdot\ell=-(1-a')S'\cdot\ell.
\end{equation}
\setcounter{subsubsection}{\value{equation}}\noindent Clearly,
$-(K_{S'}+\Diff{S'}{\De'})\cdot \ell$ depends only on
$(S',\Diff{S'}{\De'})$, but not on $\widehat X'$. So we assume
that $-(K_{S'}+\Diff{S'}{\De'})\cdot \ell$ is fixed. Recall that
the coefficients of $\Diff{S'}{\De'}$ are standard
(see~\cite[3.9]{Sh}, \cite[16.6]{Ut}), so we can write
$\Diff{S'}{\De'}=\sum_{i=1}^r (1-1/m_i)\Xi_i'$, where $m_i\in\NN$,
$r\ge 0$. Put $m':=\mt{l.c.m.}(m_1,\dots , m_r)$.
By~\cite[3.9]{Sh} both $m'S'$ and $m'(K_{S'}+\Diff{S'}{\De'})$ are
Cartier along $\ell$. So \eref{eq-diff-ell-de} can be rewritten as
$N=(1-a')k$, where $N=-m'\ell \cdot (K_{S'}+\Diff{S'}{\De'})$ is a
fixed natural number and $k=-m'(\ell\cdot S')$ is also natural.
Thus by \eref{eq-a-final} $N=(1-a')k>kr\ep\ge r\ep$. This gives us
that $r<N/\ep$ is bounded and proves the theorem.
\end{proof}

Now, we present a few corollaries of Theorem~\ref{result-Fano-1}
and Theorem~\ref{theorem-cover}. We concentrate our attention on
the three-dimensional case (then all required conjectures are
known to be true, see \cite{Sh1} and \cite{A}). Recall in this
case a non-exceptional contraction such as in
Theorem~\ref{result-Fano-1} has either $1$, $2$, $3$, $4$, or
$6$-complement.

Put $X=Z$ and $D=0$ in Theorem~\ref{theorem-cover}. We obtain
\begin{corollary}
\label{singularity-pi-1}
Let $(Z\ni o,D)$ be a three-dimensional exceptional klt
singularity such that $\totaldiscr{Z,D}>-1+\ep$ and $D\in\Msm$.
Then
\para
\label{singularity-pi-1-i}
the order of algebraic fundamental group
$\pi_1^{\mt{alg}}(Z\setminus \Sing(Z))$ is bounded by a constant
$\Const(\ep)$;
\para
the index of $K_Z+D$ is bounded by a constant $\Const(\ep)$;
\para
for any exceptional divisor $E$ over $Z$ we have either
$\dis{E}>0$ or $\dis{E}\in \mathfrak{M}(\ep)$, where
$\mathfrak{M}(\ep)\subset (-1,0]$ is a subset which depends only
on $\ep$.
\end{corollary}
Note that without assumption of exceptionality,
$\pi_1^{\mt{alg}}(Z\setminus \Sing(Z))$ is not bounded, however it
is finite \cite[Th. 3.6]{SW}. The assertion of
\eref{singularity-pi-1-i} also holds for the topological fundamental
group $\pi_1$ under the assumption that $\pi_1(Z\setminus \Sing(Z))$
is finite. M.~Reid has informed us that the finiteness of
$\pi_1(Z\setminus \Sing(Z))$ for three-dimensional log terminal
singularities was proved by N.~Shepherd-Barron (unpublished).

\begin{corollary}[\cite{Pr2}]
Fix $\ep>0$. Let $(X/Z\ni o,D)$ be a three-dimensional log variety
of local type such that $K+D$ is $\QQ$-Cartier and $-(K+D)$ is
$f$-nef and $f$-big. Assume that $f$ is exceptional and
$\totaldiscr{X}>-1+\ep$.
\para
If $\dim(Z)\ge 2$, then $\pi_1^{\mt{alg}}(Z\setminus \Sing(Z))$ is
bounded by a constant $\Const(\ep)$.
\para
If $\dim(Z)=1$, then the multiplicity of the central fiber
$f^{-1}(o)$ is bounded by a constant $\Const(\ep)$.
\end{corollary}

\begin{corollary}[\cite{Sh1}]
\label{cor-small}
Fix $\ep>0$. Let $(X/Z\ni o,D)$ be a three-dimensional exceptional
log pair such that the structure morphism $f\colon X\to Z\ni o$ is
a small contraction (i.e. $f$ contracts only a finite number of
curves), $\totaldiscr{X,D}>-1+\ep$, $D\in \Mm^3$ and $-(K+D)$ is
nef and big over $Z$. Then
\para
\label{corollary-rho-i}
$\rho(X/Z)$ and $\rho^{\mt{an}}(X/Z)$ are bounded by
$\Const(\ep)$;
\para
\label{corollary-rho-ii}
the number of components of the central fiber $f^{-1}(o)$ is
bounded by $\Const'(\ep)$.
\end{corollary}
\begin{proof}
Notation as in the proof of Theorem~\ref{result-Fano-1}. Take some
$n$-complement $K_{\widehat X}+\De+S+\Upsilon$ with $n\le N_2$.
Run $(K_{\widehat X}+\De+\Upsilon)$-MMP. For each extremal ray $R$
we have $R\cdot S>0$. Hence $S$ is not contracted. At the end we
get a model $p\colon \widetilde X\to Z$ with $p$-nef
$K_{\widetilde X}+\widetilde \De+\widetilde\Upsilon\equiv
-\widetilde S$. Since $K_{\widehat X}+\De+S+\Upsilon$ is
numerically trivial, for any divisor $E$ of $\KKK(X)$, we have
$\dis{E,\De+S+\Upsilon}=\dis{E,\widetilde \De+\widetilde
S+\widetilde\Upsilon}$ (cf. \cite[3.10]{Ko}). This shows that
$K_{\widetilde X}+\widetilde \De+\widetilde S+\widetilde\Upsilon$
is plt. Further, by Lemma~\ref{dim-q(S)>0}, $p(\widetilde S)=o$.
Since $-\widetilde S$ is nef over $Z$, we see that $\widetilde S$
coincides with the fiber over $o$. By construction,
$n\pair{K_{\widetilde S}+\Diff{\widetilde S}{\widetilde
\De+\widetilde\Upsilon}}\sim 0$ and $K_{\widetilde
S}+\Diff{\widetilde S}{\widetilde \De+\widetilde\Upsilon}$ is klt
(by the Adjunction \cite[17.6]{Ut}). Therefore
\[
\totaldiscr{\widetilde S, \Diff{\widetilde S}{\widetilde
\De+\widetilde\Upsilon}}\ge -1+1/n, \quad n\le N_2.
\]
Obviously, $\Diff{\widetilde S}{\widetilde
\De+\widetilde\Upsilon}\ne 0$. By \cite{A}, $\widetilde S$ belongs
to a finite number of algebraic families. Thus we may assume that
$\rho(\widetilde S)$ is bounded by $\Const(\ep)$.
\par
Now, consider the exact sequence
\[
0\longrightarrow\ZZ\longrightarrow \OOO^{\mt{an}}_{\widetilde X}
\stackrel{\exp}{\longrightarrow} \OOO^{\mt{an}*}_{\widetilde
X}\longrightarrow 0.
\]
By Kawamata-Viehweg vanishing $R^if^*\OOO^{\mt{an}}_{\widetilde
X}=0$, $i>0$. Hence, $\Pic^{\mt{an}}\left({\widetilde
X}\right)=H^2\left({\widetilde X},\ZZ\right)$. Similarly,
$H^2\left(\widetilde S,\ZZ\right)=\Pic\left(\widetilde S\right)$.
Since ${\widetilde X}$ is a topological retract of $\widetilde
S=p^{-1}(o)$, $H^2({\widetilde X},\ZZ)=H^2(\widetilde S,\ZZ)$.
Hence $\rho^{\mt{an}}(\widetilde X)$ is bounded, and so is
$\rho^{\mt{an}}(\widehat X)$ (because $\widehat X\bir \widetilde
X$ is a sequence of flips). This shows \eref{corollary-rho-i}. To
prove \eref{corollary-rho-ii} one can use that
$\rho^{\mt{an}}(X/Z)$ is equal to the number of components of
$f^{-1}(o)$ (by the same arguments as above, see
\cite[(1.3)]{Mo}).
\end{proof}

\begin{corollary}
Fix $\ep>0$. Let $(Z\ni o,D)$ be a three-dimensional exceptional
klt singularity such that $\totaldiscr{X,D}>-1+\ep$ and $D\in
\Mm^3$. Then for its $\QQ$-factorialization $f\colon X\to Z$ one
has
\para
$\rho(X/Z)$ and $\rho^{\mt{an}}(X/Z)$ are bounded by
$\Const(\ep)$;
\para
the number of components of $f^{-1}(o)$ is bounded by
$\Const'(\ep)$.
\end{corollary}

Note that for non-exceptional flipping contractions the number of
components of the fiber is not bounded even in the terminal case
\cite[13.7]{KoM}. We present an example of a flopping contraction as
in Corollary~\ref{cor-small}:

\begin{example}
Let $(Z\ni o)$ be a hypersurface singularity given by
$x_1^3+x_2^3+x_3^5+x_4^5$ in $\CC^4$. By \cite{IP} $(Z\ni o)$ is
exceptional (and canonical). It is easy to see also that it is not
$\QQ$-factorial. Let $f\colon X\to Z$ be a $\QQ$-factorialization
\cite[6.11.1]{Ut}. By Lemma~\ref{lemma-crepant-ex}, $(X/Z\ni o,0)$
is exceptional. Hence it satisfies conditions of
Corollary~\ref{cor-small} (with $D=0$).
\end{example}

Many examples of exceptional singularities can be found in
\cite{MP} and \cite{IP}. Finally, we propose an example of an
exceptional Fano contraction $f\colon X\to Z$ with
$\dim(X)>\dim(Z)$.

\begin{example}[{\cite[Sect. 7]{Pr3}}]
Starting with $\PP^1\times\CC^1$, blow-up points on a fiber of the
projection $\PP^1\times\CC^1\to \CC^1$ so that we obtain a
fibration $f^{\min}\colon X^{\min}\to \CC^1$ with the central
fiber having the following dual graph
\[
\begin{array}{ccccccccccc}
 &&&&\stackrel{-3}{\circ}&&&&&&\\
 &&&&|&&&&&&\\
\stackrel{-2}{\circ}&\pal&\stackrel{-2}{\circ}&\pal&\stackrel{-b}{\circ}
&\pal&\stackrel{-2}{\circ}&\pal&\stackrel{-1}{\bullet}&\pal&
\underbrace{\stackrel{-3}{\circ} \pal\stackrel{-2}{\circ}\pal
\cdots\pal\stackrel{-2}{\circ}}_{b-2}\\
\end{array}
\]
where $b\ge 2$. Now, contract curves corresponding to white
vertices. We obtain an extremal contraction $f\colon X\to \CC^1$
with two log terminal points. The canonical divisor $K_X$ is
$3$-complementary, but not $1$ or $2$-complementary \cite[Sect.
7]{Pr3}. Hence $f$ is exceptional.
\end{example}

\section{Global case}
In this section we modify Theorem~\ref{result-Fano-1} to the
global case. In contrast with the local case here we have to
assume also the existence of a boundary with rather ``bad''
singularities. Theorem~\ref{result-Fano-2} is a special case of
Conjecture~\ref{conjecture-inductive-complements}.
\begin{theorem}[Global case]
\label{result-Fano-2}
Let $(X,D)$ be a $d$-dimensional log variety of global type such
that
\para
$K+D$ is klt;
\para
$-(K+D)$ is nef and big;
\para
$D\in\MMM$, where $\MMM=\Mm^d$ or $\Msm$.
\par \smallskip
Assume that there is a boundary $D^\flat$ such that
\para
$K+D+D^\flat$ is not klt;
\para
$-(K+D+D^\flat)$ is nef and big.
\par\smallskip
Assume LogMMP in dimension $d$. Then there exists a non-klt
$n$-complement of $K+D$ for $n\in\NNN_{d-1}(\MMM)$.
\end{theorem}
\begin{proof}
First, replace $X$ with its $\QQ$-factorialization. Then as in
Lemma~\ref{weak-log-Fano} we take $D^\mho\ge 0$ such that
$-(K+D+D^\flat+D^\mho)$ is ample (but $K+D^\flat+D^{\flat\mho}$ is
not necessarily lc). Next we put
$D^\triangledown=t(D^\flat+D^\mho)$, $0<t\le 1$ so that
\para
$K+D+D^\flat$ is lc but not klt (i.e. $t$ is the log canonical
threshold $c(X,D,D^\flat+D^{\mho})$).
\par\smallskip
Now, the proof is similar to the proof of
Theorem~\ref{result-Fano-1}.
\end{proof}

\begin{corollary}[cf. {\cite[2.8]{Sh1}}]
Let $(X,D)$ be a $d$-dimensional log variety of global type such
that
\para
$K+D$ is klt;
\para
$-(K+D)$ is nef and big;
\para
$D\in\MMM$, where $\MMM=\Mm^d$ or $\Msm$.
\para
$(K+D)^d>d^d$.
\par\smallskip
Assume LogMMP in dimension $d$. Then there exists a non-klt
$n$-complement of $K+D$ for $n\in\NNN_{d-1}(\MMM)$.
\end{corollary}
\begin{proof}
A boundary $D^\flat$ such as in Theorem~\ref{result-Fano-2} exists
by Riemann-Roch (see e.~g. \cite[6.7.1]{Ko}).
\end{proof}

Many examples of exceptional log del Pezzo surfaces can be found
in \cite{Sh1}, \cite{Abe}, \cite{KeM} and \cite{Pr3}.

\section{Appendix}
In this section we give two very useful properties of complements.
We will use Definition \eref{def-coplements-n} for the case when
$D$ is a subboundary, i.e. a $\QQ$-divisor (not necessarily
effective) with coefficients $d_i\le 1$.
\begin{proposition1}[{\cite[2.13]{Sh1}}]
\label{bir-prop}
Fix $n\in\NN$. Let $f\colon Y\to X$ be a birational contraction
and let $D$ be a subboundary on $Y$ such that
\para
\label{prop-i}
$K_Y+D$ is nef over $X$;
\para
\label{prop-ii}
$f(D)=\sum d_if(D_i)$ is a boundary whose coefficients satisfy the
inequality
\[
\down{(n+1)d_i}\ge nd_i.
\]
Assume that $K_X+f(D)$ is $n$-complementary. Then $K_Y+D$ is also
$n$-complementary.
\end{proposition1}

\begin{proof}
Let us consider the crepant pull-back $K_Y+D'=f^*(K_X+f(D)^+)$,
$f_*D'=D$. Write $D'=S'+B'$, where $S'$ is reduced, $S'$, $B'$
have no common components, and $\down{B'}\le 0$. We claim that
$K_Y+D'$ is an $n$-complement of $K_Y+D$. The only thing we need
to check is that $nB'\ge \down{(n+1)\fr{D}}$. From \eref{prop-ii}
we have $f(D)^+\ge f(D)$. This gives us that $D'\ge D$ (because
$D-D'$ is $f$-nef; see \cite[1.1]{Sh}). Finally, since $nD'$ is an
integral divisor, we have
\[
nD'\ge nS'+\down{(n+1)B'}\ge n\down{D}+\down{(n+1)\fr{D}}.
\]
\end{proof}

The following is a refinement of \cite[Proof of 5.6]{Sh} and
\cite[19.6]{Ut}.

\begin{proposition1}[{\cite{Pr2}}]
\label{prodolj}
Fix $n\in\NN$. Let $(X/Z\ni o,D=S+B)$ be a log variety. Set
$S:=\down{D}$ and $B:=\fr{D}$. Assume that
\para
$K_X+D$ is plt;
\para
$-(K_X+D)$ is nef and big over $Z$;
\para
$S\ne 0$ near $f^{-1}(o)$;
\para
\label{prop-2-iv}
the coefficients of $D=\sum d_iD_i$ satisfy the inequality
\[
\down{(n+1)d_i}\ge nd_i.
\]
Further, assume that near $f^{-1}(o)\cap S$ there exists an
$n$-complement $K_S+\Diff SB^+$ of $K_S+\Diff SB$. Then near
$f^{-1}(o)$ there exists an $n$-complement $K_X+S+B^+$ of
$K_X+S+B$ such that $\Diff SB^+=\Diff S{B^+}$.
\end{proposition1}

\begin{proof}
Let $g\colon Y\to X$ be a log resolution. Write
$K_Y+S_Y+A=g^*(K_X+S+B)$, where $S_Y$ is the proper transform of
$S$ on $Y$ and $\down{A}\le 0$. By the Inversion of Adjunction
\cite[17.6]{Ut}, $S$ is normal and $K_S+\Diff SB$ is plt. In
particular, $g_S\colon S_Y\to S$ is a birational contraction.
Therefore we have
\[
K_{S_Y}+\Diff{S_Y}{A}=g_S^*(K_S+\Diff SB).
\]
Note that $\Diff{S_Y}{A}=A|_{S_Y}$, because $Y$ is smooth. It is
easy to show (see \cite[4.7]{Pr3}) that the coefficients of $\Diff
SB$ satisfy the inequality \eref{prop-2-iv}. So we can apply
Proposition~\ref{bir-prop} to $g_S$. We get an $n$-complement
$K_{S_Y}+\Diff{S_Y}A^+$ of $K_{S_Y}+\Diff{S_Y}A$. In particular,
by \eref{def-coplements-n}, there exists
\[
\Theta\in \left|-nK_{S_Y}-\down{(n+1) \Diff{S_Y}A}\right|
\]
such that
\[
n\Diff{S_Y}A^+= \down{(n+1)\Diff{S_Y}A}+\Theta.
\]
By Kawamata-Viehweg Vanishing,
\begin{multline*}
R^1h_*(\OOO_{Y}(Y,-nK_Y-(n+1)S_Y-\down{(n+1)A}))=\\
R^1h_*(\OOO_{Y}(Y, K_Y+\up{-(n+1)(K_Y+S_Y+A)}))=0.
\end{multline*}
 From the exact sequence
\begin{multline*}
0\longrightarrow\OOO_{Y}(Y,-nK_Y-(n+1)S_Y-\down{(n+1)A})\\
\longrightarrow\OOO_{Y}(Y,-nK_Y-nS_Y-\down{(n+1)A})\\
\longrightarrow\OOO_{S}(S_Y,-nK_{S_Y}-\down{(n+1)A}|_{S_Y})
\longrightarrow 0
\end{multline*}
we get surjectivity of the restriction map
\begin{multline*}
H^0(Y,\OOO_{Y}(-nK_Y-nS_Y-\down{(n+1)A})) \longrightarrow\\
H^0(S_Y,\OOO_{S_Y}(-nK_{S_Y}-\down{(n+1)A}|_{S_Y})).
\end{multline*}
Therefore there exists a divisor
\[
\Xi\in\left|-nK_Y-nS_Y-\down{(n+1)A}\right|
\]
such that $\Xi|_{S_Y}=\Theta$. Set
\[
A^+:=\frac1{n}(\down{(n+1)A}+\Xi).
\]
Then $n(K_Y+S_Y+A^+)\sim 0$ and $(K_Y+S_Y+A^+)|_{S_Y}=
K_{S_Y}+\Diff{S_Y}A^+$. Note that we cannot apply the Inversion of
Adjunction on $Y$ because $A^+$ can have negative coefficients. So
we put $B^+:=g_*A^+$. Again we have $n(K_X+S+B^+)\sim 0$ and
$(K_X+S+B^+)|_S=K_S+\Diff SB^+$. We have to show only that
$K_X+S+B^+$ is lc. Assume that $K_X+S+B^+$ is not lc. Then
$K_X+S+B+\alpha(B^+-B)$ is also not lc for some $\alpha<1$. It is
clear that $-(K_X+S+B+\alpha(B^+-B))$ is nef and big over $Z$. By
the Inversion of Adjunction \cite[17.6]{Ut}, $K_X+S+B+\alpha(B^+-B)$
is plt near $S\cap f^{-1}(o)$. Hence $LCS(X,B+\alpha(B^+-B))=S$
near $S\cap f^{-1}(o)$. On the other hand, by the Connectedness Lemma
\cite[17.4]{Ut}, $LCS(X,B+\alpha(B^+-B))$ is connected near
$f^{-1}(o)$. Thus $K_X+S+B+\alpha(B^+-B)$ is plt. This
contradiction proves the proposition.
\end{proof}

\end{document}